\RequirePackage{fix-cm}
\documentclass[smallextended]{svjour3}       % onecolumn (second format)
\smartqed  % flush right qed marks, e.g. at end of proof
\usepackage{graphicx}
\usepackage[margin=1in]{geometry} 
\usepackage{amsmath,amssymb}
\usepackage{hyperref}
\usepackage{caption}
\usepackage{soul}
\usepackage{tcolorbox}
\usepackage{float}
\definecolor{red}{rgb}{1,0,0} 
\newcommand{\bV}{\mathbf{V}}
\newcommand{\bv}{\mathbf{v}}
\newcommand{\bu}{\mathbf{u}}
\newcommand{\bn}{\mathbf{n}}
\newcommand{\bb}{\mathbf{b}}
\newcommand{\by}{\mathbf{y}}
\newcommand{\ignore}[1]{}

\ignore{
\newtheorem{theorem}{Theorem}

\newtheorem{definition}{Definition}

\newtheorem{proposition}{Proposition}
}

\begin{document}

\title{A Note on the Stability of  Two Families of Two-Step Schemes%\thanks{Grants or other notes
}


\author{Xiaoming Wang         \and
        Yinqian Yu %etc.
}


\institute{Xiaoming Wang \at
             Department of Mathematics and Statistics, Missouri University of Science and Technology, Rolla, MO 65409, USA\\
             % National Center for Applied Mathematics Shenzhen (NaCAMS), ShenZhen 518055, PR China\\
             % Guangdong Provincial Key Laboratory of Computational Science and Material Design, Southern University of Science and Technology, ShenZhen 518055, PR China \\
             % SUSTech International Center for Mathematics, Southern University of Science and Technology, ShenZhen 518055, PR China \\
             % Tel.: +123-45-678910\\
             % Fax: +123-45-678910\\
              \email{wxm.math@outlook.com}           %  \\
           \and
           Yinqian Yu \at
              Department of Mathematics, The University of Utah, Salt Lake City, UT 84112, USA
}

\date{Received: date / Accepted: date}

\maketitle

% Uncomment and use as if needed
%\newtheorem{theorem}{Theorem}
%\newtheorem{lemma}[theorem]{Lemma}
%\newdefinition{rmk}{Remark}
%\newproof{pf}{Proof}
%\newproof{pot}{Proof of Theorem \ref{thm}}

%\cortext[1]{Corresponding author, orcid=0000-0002-2399-6336}

%\fntext[1]{1}


\begin{abstract}
We investigate the stability of two families of three-level two-step schemes that extend the classical second order BDF (BDF2) and second order Adams-Moulton (AM2) schemes. For a free parameter restricted to an appropriate range that covers the classical case, we show that both the generalized BDF2 and the generalized AM2 schemes are A-stable. 
We also introduce the concept of  uniform-in-time stability which characterizes a scheme's ability to inherit the uniform boundedness over all time of the solution of  damped and forced equation with the force uniformly bounded in time. We then demonstrate that 
 A-stability and uniform-in-time stability are equivalent for  three-level two-step schemes. %when applied to damped and forced equation with the force uniformly bounded in time.
Next, these two families of schemes are utilized to construct efficient and unconditionally stable IMEX schemes for systems that involve a damping term, a skew symmetric term, and a forcing term. These novel IMEX schemes are shown to be uniform-in-time energy stable in the sense that the norm of any numerical solution is bounded uniformly over all time, provided that the forcing term is uniformly bounded time, the skew symmetric term is dominated by the dissipative term, together with a mild time-step restriction.
% The IMEX algorithms are applied to the Stokes-Darcy system leading to a new family of uniform in time energy stable scheme.
Numerical experiments verify our theoretical results. They also indicate that the generalized schemes could be more accurate  and/or  more stable than the classical ones for suitable choice of the parameter.
\end{abstract}



% Keywords

\maketitle

% Main text
\section{Introduction}
\label{intro}
Efficiency of numerical simulation calls for higher order numerical schemes which enjoy smaller local truncation error for the same step-size when compared to lower order methods. This allows, at least heuristically, relatively larger time-step without committing larger local truncation error, and hence one can reach the terminal time faster under the same error tolerance. 
There are two popular types of higher order in time numerical methods: the Runge-Kutta method and the multi-step method. We will focus on the multi-step method in this manuscript.

The main idea of multi-step method is to utilize numerical approximations of the exact solution at multiple previous times to construct approximate solution at a future time close to the present one.
A general linear multi-step method (LMM) with $s$-steps applied to $y'=f(t,y)$ takes the form
% s-step methos
\begin{equation}\label{general_sstep_method}
    \sum_{m=0}^s a_m y_{n+m} = h \sum_{m=0}^sb_mf(t_{n+m},y_{n+m}),
\end{equation}
with $a_s=1$ (or $a_s>0$), where $h$ is the uniform time step size, $t_n=nh$ assuming the starting time is $t=0$, and $y_n$ represents the numerical approximation of $y(nh)$\cite{wanner1996solving}.
The method is characterized in terms of the two (generating) polynomials:
% 1st and 2nd char poly
\begin{equation} \label{gen-poly}
    \rho(w) = \sum_{m=0}^s a_mw^m, \quad \sigma(w) = \sum_{m=0}^sb_mw^m.
\end{equation}
It is easy to imagine that an $s$-step method that involves $s+1$-levels of approximate solutions can be of the order of $s$ according to Taylor's theorem. 
We recover the Backward Differentiation Formula (BDF) method if the expansion is centered about the future time grid together with  $b_s=1$. We arrive at the classical second Adams-Moulton (AM) method if the mid-point of the current and future time grid is utilized for the expansion and $b_s=b_{s-1}=\frac12$ in the second order case \cite{hairer1993solving,iserles2009first,wanner1996solving}.

One can extend the two numerical schemes by allowing more general combinations of the terms on the right hand side without sacrificing the accuracy,  and this is how we generalize the classical BDF2 and AM2 schemes.  
While the extensions presented in this manuscript maintain the two-step (three-level) structure, more general approaches that involve additional levels in developing new algorithms by allowing combinations different from the classical ones have been utilized in other works in various context, see for instance \cite{chen2013efficient,cheng2016long,gottlieb2012long,gottlieb2012stability,shen2023split}.

Besides the local truncation error consideration, 
the stability of numerical schemes is of great importance in the performance of the algorithm. 
Dahlquist's equivalence theorem states that the convergence of a numerical scheme is equivalent to the consistency plus the stability of the algorithm.
Indeed, it is easy to construct numerical algorithms that have nice local truncation error but fail to well-approximate the exact solution
 \cite{hairer1993solving,iserles2009first,wanner1996solving}.
A variety of stability concepts have been proposed and investigated. The classical treatise by Hairer and Wanner offers a comprehensive survey of related topics \cite{wanner1996solving}.
One of the most important stability concepts is the so-called A-stability, which requires the stability region of the numerical scheme cover the left half-plane when applied to the scalar ODE $y'=\lambda y$ with the real part of $\lambda$ non-positive. 
It turns out that the generalized BDF2 and AM2 schemes developed by including a free parameter  retain the A-stability when the free parameter falls into an appropriate regime.
Our numerics also indicate that the generalized schemes could be more stable or more accurate than the classical ones for suitable value of the parameter.

On the other hand, for applications to dissipative ODEs or PDEs, the stability is often in terms of the boundedness of certain norms of the numerical solution. In addition, when the physical phenomena involve long-time behavior such as the climate, the coarsening process associated with some material science problems, and transport in porous media, it is natural to investigate long-time stability. Long-time stability issue of numerical schemes has been investigated by many authors in different contexts, see for instance
\cite{chen2013efficient,chen2016efficient,cockburn1997estimating,foias1989exponential,gottlieb2012long,shen2023split,shen1990long,shen2018scalar,stuart1998dynamical,stuart1994numerical,yang2017numerical}. 
This motivates us to introduce the concept of uniform-in-time (energy) stability of numerical schemes, i.e., whether the numerical method can inherit the uniform boundedness of the solution over all time of the underlying linear damped-forced ODE $y'=\lambda y +g$ with $\it{Re}(\lambda) < 0$ and $g\in L^\infty [0,\infty)$. 
One natural follow-up question then is the relationship between A-stability and uniform-in-time (energy) stability.

The main results of this note are
 \begin{enumerate}
\item The introduction of two families of second order in time two-step methods that are extensions of the classical backward differentiation method (BDF2) and Adams-Moulton (AM2) methods respectively; 
\item The identification of the parameter regimes that ensure A-stability of the generalized BDF2 and AM2 schemes;
\item The equivalence of uniform-in-time (energy) stability and  A-stability for three-level two-step schemes;
\item The development and analysis of efficient and stable methods, that inherit the boundedness over all time of the solutions of the underlying model,  by combining the generalized BDF2/AM2 with Implicit-Explicit (IMEX) methodology  when applied to problems in the presence of a skew symmetric term in addition to damping and forcing.
\end{enumerate}
Applications to some prototype linear dissipative problems, the Stokes-Darcy system in particular,  will be briefly discussed in the appendix.

The rest of the note is organized as follows. In the next section, we present the novel generalized BDF2 and generalized AM2 schemes. In the third section, we study the A-stability of the generalized BDF2 method and the generalized AM2.
The relationship between A-stability and uniform-in-time energy stability is analyzed in the fourth section.
Stable and efficient numerical schemes based on the generalized BDF2/AM2 and IMEX method that are applicable to models with a skew-symmetric term are presented in the fifth section.
Numerical results are reported in the sixth section.
We offer our conclusion and perspectives at the end.
Application to the Stokes-Darcy system is outlined in the appendix.

%%%%%
\section{Two families of second order three-level two-step methods}
\label{Sec2}
%\textcolor{red}{In this section, some classical and corresponding generalized second order two step methods are recalled.}
The main purpose of this section is to present two families of three-level two-step LMM methods. They are generalizations of the classical BDF2 and AM2 methods respectively.

% 2.1
\subsection{Generalized second order BDF2 schemes}
\label{Sec2_sub1}
Recall that the classical BDF2 scheme applied to ODE $y'=f(t,y)$ takes the following form.
%% BDF2
\begin{equation}\label{BDF2}
        \frac{3}{2}y_{n+1}-2y_n+\frac{1}{2}y_{n-1} = hf_{n+1},
\end{equation}
where $h$ is the uniform time-step size, and $y_n$ is the approximation of $y(t_n)$.
BDF methods are widely used for integration of stiff differential equations due to their strong stability \cite{hairer1993solving,iserles2009first}.

Notice that we could include terms involving $f_n, f_{n-1}$ on the right hand side without sacrificing the formal accuracy as long as the combined effect of the extra terms are of higher order.  Hence, we propose the following generalized BDF2 scheme
%% g-BDF2
\begin{equation}\label{generalized_BDF2_general_form}
    \frac{3}{2}y_{n+1} - 2y_n + \frac{1}{2}y_{n-1} = h (\alpha f(t_{n+1},y_{n+1}) + (2-2\alpha)f(t_n,y_n) + (\alpha-1)f(t_{n-1},y_{n-1}) ),
\end{equation}
where  $\alpha$ is a free parameter.
Notice that we recover the classical BDF2 scheme \eqref{BDF2} when $\alpha=1$.

It is easy to verify the second order accuracy of the novel generalized BDF2 scheme utilizing classical results with the aid of the generating polynomials $\rho$ and $\sigma$ for this scheme
\begin{equation}
    \rho(w) = \frac32 w^2 -2w +\frac12,  \quad \sigma(w) = \alpha w^2 + (2-2\alpha)w + (\alpha-1)
\end{equation}
It is straightforward to verify that  $\rho(w) - \sigma(w)lnw = c(w-1)^{3} + \mathcal{O}(|w-1|^{4}), \quad \mbox{as}\ w \to 1$. Hence, we have the following proposition (by Theorem 2.1 in \cite{iserles2009first}).
\begin{proposition}
The generalized BDF2 scheme \eqref{generalized_BDF2_general_form} is of order 2 for all $\alpha$.
\end{proposition}

% 2.2
\subsection{Generalized second order Adams-Moulton method}
\label{Sec2_sub2}
Recall the classical second order Adams-Moulton method (AM2) applied to ODE $y'=f(t,y)$ takes the form 
$$   y_{n+1} = y_n + h(\frac{1}{2}f_n + \frac{1}{2}f_{n+1}).$$
There is also an alternative Adams-Bashforth approach (AB2) that treats the right hand side explicitly
$$  y_{n+1} = y_n + h(\frac{3}{2}  f_n - \frac{1}{2}f_{n-1}). $$
Combining the classical (implicit) AM2 and explicit AB2 schemes, we arrive at the generalized second order Adams-Moulton method
% new scheme here
\begin{equation}\label{generalized_AM2_general_form}
   y_{n+1} -  y_n = h (\alpha f(t_{n+1},y_{n+1}) + (\frac{3}{2}-2\alpha)f(t_n,y_n) + (\alpha-\frac{1}{2})f(t_{n-1},y_{n-1}) ),
\end{equation}
where $\alpha$ is a parameter that allows us to recover the classical AM2 scheme when $\alpha=1/2$, and the explicit AB2 scheme when $\alpha=0$.

The second order accuracy of the generalized AM2 scheme can be verified in a straightforward manner as for the generalized BDF2 scheme by utilizing the generating polynomials.
\begin{proposition}
The generalized AM2 scheme \eqref{generalized_AM2_general_form} is of second order for all $\alpha$.
\end{proposition}

%%%%%
\section{A-stability of the generalized BDF2 and AM2 schemes}
\label{Sec3}
%\textcolor{red}{The uniform-in-time energy stability will be treated separately using a unified approach that is suitable for A-stable schemes.}
The goal of this section is to show that the generalized BDF2 schemes are A-stable provided that  $\alpha\ge 3/4$,  while the generalized AM2 schemes are A-stable provided $\alpha\ge 1/2$. 
%We also comment on how the stability region depends on the choice of the free parameter $\alpha$ based on our numerical experiments.

%
\subsection{A-stability of the generalized BDF2 method}
\label{Sec3_sub1}
%\textcolor{red}
{In this subsection, we identify the parameter values that guarantee the A-stability of the generalized BDF2 scheme.}

The main result of this subsection is the following theorem.
%% A-stbaility of g-BDF2
\begin{theorem}
The generalized BDF2 method \eqref{generalized_BDF2_general_form} is A-stable if and only if $ \alpha \ge \frac{3}{4}$.
\end{theorem}

\begin{proof}
The proof relies on an equivalent condition on A-stability formulated using the characteristic polynomial $\eta(z,w)=\rho(z,w)-z\sigma(z,w)$ where $\rho$ and $\sigma$ are the two generating polynomials defined in \eqref{gen-poly}. 
 
 According to Lemma 4.8 in \cite{iserles2009first} a general LMM scheme is A-stable if and only if 
\begin{equation}
    b_s > 0, \quad |w_1(it)|, |w_2(it)|,...,|w_{q(it)}(it)| \le 1, \forall t \in \mathbb{R},
\end{equation}
where $w_1,...,w_{q(z)}$\footnote{We have emphasized the dependence on $z$ of the roots of the characteristic polynomial here.} are all the zeros of the characteristic polynomial $\eta(z,w) = \rho(z,w)-z\sigma(z,w)$.

The characteristic polynomial of the generalized BDF2 method is
\begin{equation}
    \eta(z,w) = (\frac{3}{2} - \alpha z)w^2 +(-2+(2\alpha -2)z)w + \frac{1}{2}+(1-\alpha)z.
\end{equation}

%Let $z = it, t\in \mathbb{R}$, 
Hence, the A-stability of \eqref{generalized_BDF2_general_form} is equivalent to both zeros of the following quadratic function lie in the closed unit disc in the complex plane for all $ t\in \mathbb{R}$.
\begin{equation}
    (\frac{3}{2}-\alpha it)w^2 +(-2+(2\alpha -2)it)w + \frac{1}{2} + (1-\alpha)it,\quad t\in \mathbb{R}.
\end{equation}

We utilize the Cohn-Schur criterion(\cite{iserles2009first} Lemma 4.9) to check the condition of the roots of the quadratic function above. The Cohn-Schur criterion states that the condition of both zeros of the quadratic function $aw^2 + bw + c$ residing in the closed complex unit disc is equivalent to:
\begin{equation}
     |a| \ge |c|, \quad ||a|^2-|c|^2|\ge |a\Bar{b}-b\Bar{c}|.
     \label{Cohn-Schur}
\end{equation}

Notice that for the generalized BDF2 we have
\begin{equation}
    a = \frac{3}{2}-\alpha it,\quad b = -2+(2\alpha -2)it,\quad c =\frac{1}{2} + (1-\alpha)it.
\end{equation}

Hence, the equivalent condition for A-stability is 
\begin{equation}
\begin{split}
& (1-2\alpha)t^2 \le 2, \quad \forall
 t \in \mathbb{R},\\
 &  |2 + (2\alpha-1)t^2| \ge |(2-2\alpha)t^2-2+2it|,\quad \forall t \in \mathbb{R}, 
\end{split}
\end{equation}
which holds iff $\alpha \ge \frac{3}{4}$.

This end the proof of the A-stability of the generalized BDF2 scheme.
\end{proof}

%

% S3.2
\subsection{A-stability of the generalized AM2 scheme}
\label{Sec3_sub2}

%\textcolor{red}{Dido the same thing as for the BDF2 schemes.}

In this subsection, we show that the generalized AM2 scheme is A-stable if the free parameter $\alpha$ is at least one-half.

%\subsection{A-stability of the generalized AM2 method}
\begin{theorem}
The generalized AM2 scheme \eqref{generalized_AM2_general_form} is A-stable if and only if $\alpha \ge 1/2$.

\end{theorem}

\begin{proof}
We follow the same steps as the proof for the generalized BDF2 method. 

Consider the characteristic polynomial of the generalized AM2 method \eqref{generalized_AM2_general_form}
\begin{equation}
    \eta(z,w) = (1-\alpha z)w^2 + (-1 + (2\alpha -\frac{3}{2})z)w + (\frac{1}{2} - \alpha)z. \quad z = it, t\in \mathbb{R}.
\end{equation}

By Lemma 4.7 and Lemma 4.8 in \cite{iserles2009first}, the generalized AM2 scheme is A-stable if and only if all roots of $\eta(it,w), t \in \mathbb{R}$ are in the closed unit disc in the complex plane. More specifically, the scheme is A-stable if and only both zeros of the following quadratic belong to the unit disc in the complex plane.
\begin{equation}
    (1-\alpha it)w^2 + (-1 + (2\alpha -\frac{3}{2})it)w + (\frac{1}{2}-\alpha)it ,\quad t \in \mathbb{R}.
\end{equation}
 Then by the Cohn-Schur criterion(Lemma 4.9 \cite{iserles2009first}), the scheme is A-stability  if and only if \eqref{Cohn-Schur} 
%
%\begin{equation}
%     |a| \ge |c|, \quad ||a|^2-|c|^2|\ge |a\Bar{b}-b\Bar{c}|
%\end{equation}
holds for all $t \in \mathbb{R}$, where
\begin{equation}
    a = 1-\alpha it,\quad b = -1 + (2\alpha -\frac{3}{2})it, \quad c = (\frac{1}{2}-\alpha)it.
\end{equation}
We then deduce that the generalized AM2 method is A-stable if and only if
%Then we substitute these parameters $a,b,c$ into the requirement of the Cohn-Schur criterion to find 

\begin{equation}
  (\alpha -\frac{1}{4})t^2 + 1 \ge 0, \quad (\alpha-\frac{1}{2})t^2 \ge 0, \forall t\in \mathbb{R},
\end{equation}
which holds iff $\alpha \ge \frac{1}{2}$. Note that the case $\alpha = \frac{1}{2}$ is indeed the classical AM2 scheme.

This completes the proof of the A-stability of the generalized AM2 scheme when $\alpha\ge 1/2$.
\end{proof}

\noindent{\bf{Remark:}} It is of interest to investigate how the `strength` of the stability depends on the parameter. A way to quantify the `strength' is in terms of the size of the stability region. Preliminary numerical results indicate that the stability region grows as the parameter $\alpha$ increases for the generalized AM2 and the BDF2 schemes as reported in the MS dissertation of the second author \cite{yumsthesis}.The rigorous validation is not yet available so far.

%%%%%
\section{Uniform-in-time (energy) stability}
\label{Sec4}
%
%\subsection{Strict A-stability}

%
%\subsection{Strict A-stability and uniform-in-time energy stability}
%\textcolor{red}{Introduce the concept of strict A-stable schemes.
%Introduce the concept of uniform-in-time energy stability.}
While A-stability is an excellent stability concept, the implication of A-stability when applied to ubiquitous  damped-driven models in natural and engineering processes, is not clear.
For these damped and driven models, we are often interested in the long-time behavior of such systems. One prime example is the weather system whose long-time behavior is the climate. 
Therefore, it would be important to design numerical algorithms that are valid over long time. 
The development of long-time accurate numerical schemes for nonlinear models might be extremely challenging due to potential  intrinsic chaos within the system\cite{wang2010approximation,wang2016numerical} except for systems with strong hyperbolic structures. 
For systems that are dissipative in the sense that all solutions remain bounded in certain norm, we are naturally interested in designing algorithms that inherit the boundedness property of the system, even in the presence of an external forcing. This motivates the following stability concept.

%%%
\begin{definition}[Uniform-in-time (energy) stability]
A numerical method is called uniform-in-time (energy) stable if the numerical solutions remain uniformly bounded in time when applied to $y'=\lambda y + g$ with $g\in L^\infty([0,\infty)),{ \it{Re}}(\lambda)< 0$. 
Moreover, the long-time bound of the numerical solution should be asymptotically independent of the initial data at large time and small time-step.
%and asymptotically independent of the time step for small time-step sizes.
\ignore{ multi-step method \eqref{general_sstep_method} is said to be energy stable in $[0,T]$, if for such given time $T$, the numerical solution obtained from this multi-step method applied on a initial value problem is uniformly bounded in $[0,T]$. Further, if it is uniformly bounded in $[0,\infty)$, then this method is uniform-in-time energy stable. 
}
\end{definition}
\noindent{\bf Remark:} 
The essense of the definition above is the numerical scheme's ability to inherit the uniform-in-time boundedness of the solution to the underlying ODE. 

\noindent{\bf Remark:} The word `energy' is included since we usually utilize a norm to characterize the boundedness when the underlying phase space is a Hilbert space. Moreover, in the case when the unknown is the velocity field, one-half of the $L^2$ (or $l^2$) norm squared is the kinetic energy density per unit mass, and hence the  stability in norm is ``energy'' stability. 

\noindent{\bf Remark:}  Such kind of boundedness could lead to uniform-in-time error estimates for many linear systems.   See \cite{chen2013efficient,chenwb2016efficient} among others, for applications to the Stokes-Darcy system.

\noindent{\bf Remark:} We also emphasize that the {\it uniform-in-time stability} is different from the {\it unconditional stability} concept encountered in the numerical PDE circle. A scheme is called unconditional stable if all numerical solutions remain bounded on any finite time interval $[0,T]$  without any time-step restriction. The scheme is called conditionally stable if a time-step restriction is needed to ensure boundedness on some finite time interval. In contrast, the uniform-in-time stability imposes boundedness over the infinite time interval $[0,\infty)$, which is a more stringent criterion.

A natural follow-up question is how this new concept of uniform-in-time (energy) stability relates to the classical A-stability concept.
Recall that  A-stability of a general s-step method implies that all numerical solutions to the  damped unforced model $y'=\lambda y$ would converge to zero as time approaches infinity. 
In the simple case of a damped-forced model with a time-independent external force $f$,  we may easily perform a change of variable to deduce that the solution would remain uniformly bounded in time for the damped and forced problem. 
However, the equivalence in the case with a bounded time-dependent force is less straightforward.
Our main result in this section is that A-stability and uniform-in-time energy stability for three-level two-step methods are in fact equivalent.  
The main ingredient of the proof is the equivalence between the eigenvalues of the extended dynamical system on the product space induced by the LMM scheme and the roots of the scheme's characteristic polynomial.

%\textcolor{red}{Prove strict A-stable scheme is uniform-in-time energy stable. (Vector form... Energy method can be postponed.)}
%%%%% Therorem on equivalence between A-stability and long-time stability
\begin{theorem}[Equivalence between A-stability and Uniform-in-time (energy) stability]\label{Strictly_A-stability_Implies_Energy_stability}
 A convergent linear second-order three-level two-step method is A-stable if and only if it is uniform-in-time (energy) stable.
 \end{theorem}
\begin{proof}
We demonstrate the sufficiency, i.e., A-stability implies uniform-in-time (energy) stability first.

\ignore{
Let 
%%% show a general three level 2 steps method here
be a general three-level two-steps method with the first and second characteristic polynomial given by ($\rho,\sigma$), where
\begin{equation}
    \rho(w) = \alpha_2 w^2 + \alpha_1 w + \alpha_0, \quad \sigma(w) = \beta_2 w^2 +\beta_1 w+ \beta_0.
\end{equation}

Recall that the characteristic polynomial of the scheme is given by
% char poly
\begin{equation}
\eta(z,w)=\sum_{m=0}^2 (\alpha_m-\beta_m z)w^m.
\end{equation}
where $\alpha_j, \beta_j$ are the coefficients of the linear multi-step method (LMM).
}
Recall that the A-stability of a LMM scheme implies that the zeros of the characteristic polynomial $\eta(z,w)=\rho(z,w)-z\sigma(z,w)$, denoted $w_1(z), w_2(z)$,  lie within the interior of the unit disk (Lemma 4.7 in \cite{iserles2009first}).

Introducing the notation $\bV_{n+1}=(y_n, y_{n+1})^T$, the three-level two-step linear scheme can be rewritten as
\begin{equation}
    \bV_{n+1} = A\bV_n + h\bb_n,\label{extended_form}
\end{equation}
where\footnote{Note that $\beta_0g_n + \beta_1g_{n+1}+ \beta_2g_{n+2}$ is a second order approximation of $g(\beta_0 t_{n} + \beta_1 t_{n+1} + \beta_2 t_{n+2})$. Indeed, this quantity can be replaced by any other second order approximation of $g(\beta_0 t_{n} + \beta_1 t_{n+1} + \beta_2 t_{n+2})$.}
$$ A = \begin{pmatrix}
0 & 1 \\ \frac{-\alpha_0+ \beta_0 z}{\alpha_2 - \beta_2 z} & \frac{-\alpha_1+ \beta_1 z}{\alpha_2 - \beta_2 z} 
\end{pmatrix} , \quad  \bb_n = \begin{pmatrix}
0 \\ \frac{\beta_0g_n + \beta_1g_{n+1}+ \beta_2g_{n+2}}{\alpha_2 - \beta_2z}
\end{pmatrix}, \quad z =\lambda h.
$$

It is straightforward to show that the characteristic polynomial of the matrix $A$ is given by
$
    \frac{\eta(z,\Lambda)}{\alpha_2 -\beta_2 z}.
$
Hence, the characteristic polynomial of $A$ has the same roots as the characteristic polynomial of this linear multi-step method $\eta(z,\Lambda)$,  since $|\alpha_2-\beta_2 z|\ge |\alpha_2|=1, \forall z\in {\mathbb{C}}^{-}$, and since A-stability implies $\beta_2>0$ (see for instance, Lemma 4.8 of \cite{iserles2009first}).   

Recall that the roots of any A-stable LMM method lie in the interior of the unit disk in the complex plane, see for instance Lemma 4.7 of \cite{iserles2009first}.% {\color{red}cite exact result about roots}. 
Hence, the spectral radius $\rho_s$ of $A$ is less than one,  i.e., $\rho_s(A)<1$.
 This further implies, since $\rho_s(A)=\lim_{n\rightarrow\infty}\|A^n\|^{1/n}<1$,  there exists $N$, such that
 $$\|A^n\| \le \left(\frac{1+\rho_s(A)}{2}\right)^n, \forall n\ge N.$$
 Notice, simple induction implies,
 $$ \bV_{n} = A^n\bV_0 + h\sum_{j=0}^{n-1}A^{n-j-1}\bb_j.$$
 Let $b_\infty=\sup_{n\ge}\|\bb_n\|$. We observe that
 $$b_\infty \le \frac{(|\beta_2|+|\beta_1|+|\beta_0|)\|g\|_{L^\infty}}{\alpha_2}.$$
  Hence, for $n\ge N$,
 %%%
% {\color{red}To add the estimates here}

\begin{equation}
    \begin{split}
        \|\bV_n\|  \le  & \|A^n\|\|\bV_0\| + \sum_{j=0}^{n-1}\| A^{n-j-1}\|\|h\bb_j\| \\
        \le & (\frac{1+\rho_s(A)}{2})^n \|\bV_0\| + h\sum_{j=0}^{n-N-1}\|A^{n-j-1}\| b_{\infty} + h \sum_{j=n-N}^{n-1}\|A^{n-j-1}\|b_{\infty}\\
        \le &  (\frac{1+\rho_s(A)}{2})^n \|\bV_0\| + h b_{\infty} \sum_{j=0}^{n-N-1}(\frac{1+\rho(A)}{2})^{n-j-1} + h b_{\infty} \sum_{j=n-N}^{n-1}\|A\|^{n-j-1}\\
        = & (\frac{1+\rho_s(A)}{2})^n \|\bV_0\| +  h b_{\infty} (\frac{1+\rho_s(A)}{2})^N \cdot \frac{1-(\frac{1+\rho_s(A)}{2})^{n-N}}{\frac{1-\rho_s(A)}{2}} +  h b_{\infty}\sum_{l=0}^{N-1}\|A\|^l.    
    \end{split}
\end{equation}

 We then conclude that $\|\bV_n\|$ is uniformly bounded for all $n$.
 Moreover, we have
 $$\lim_{n\rightarrow\infty}\|\bV_n\| \le    h\ b_{\infty} (\frac{1+\rho_s(A)}{2})^N \cdot \frac{1}{\frac{1-\rho_s(A)}{2}} +  h\ b_{\infty} \sum_{l=0}^{N-1}\|A\|^l,$$
 which is independent of the initial data. 
 
 In order to show that $\|A\|$ is independence of  $z\in\mathbb{Z}^-$, we observe that
 \begin{eqnarray*}
   |\alpha_2-\beta_2 z|
    &\ge&\sqrt{ |\alpha_2-\beta_2\it{Re}(z)|^2+\beta_2^2\it{Im}(z)^2}
    \\
   &=&\sqrt{ |1-\beta_2\it{Re}(z)|^2+\beta_2^2\it{Im}(z)^2}
   \\
   &\ge& \max\{1, \beta_2 |z|\}, 
   \\
 |\alpha_0-\beta_0 z|&\le& |\alpha_0|+|\beta_0||z|,
 \\
    |\alpha_1-\beta_1 z|&\le& |\alpha_1|+|\beta_1||z|.
    \end{eqnarray*}
    We deduce, when combined with the fact that $\beta_2>0$ for A-stable schemes, that  $\|A\|$ is independent of $z\in \mathbb{C}^-$.
    
 This completes the proof of the sufficiency.

For the necessity, we assume that the scheme is uniform-in-time (energy) stable. This implies that the scheme, when formulated in the product space as in \eqref{extended_form}, must have all eigenvalues of the matrix $A$ contained in the interior of the unit disk in the complex plane, i.e., $\rho_s(A)<1$. Utilizing the fact that the eigenvalues of $A$ and the eigenvalues of the characteristic polynomial $\eta$ of the scheme are identical, we conclude that the scheme is A-stable thanks to Lemma 4.7 in \cite{iserles2009first}.
This completes the proof of the necessity.

\end{proof}

\noindent\textbf{Remark:} For the special case of generalized BDF2 and generalized AM2 schemes that we introduced in the previous sections, we can have explicit form of the matrix $A$, including its eigenvalues and (generalized) eigenvectors. Utilizing the detailed information, we can show that the long-time bound of the norm of the numerical solution is consistent with the long-time bound of the original model $y'=\lambda y + g$.
More details can be found in the thesis of the second author \cite{yumsthesis}.

\noindent\textbf{Remark:} The equivalence between A-stability and uniform-in-time stability stops at second order multi-step schemes due to Dahlquist's second barrier, which dictates that the highest order possible of an A-stable multi-step scheme is 2. On the other hand, we can infer the existence of uniform-in-time energy stable third-order multi-step schemes \cite{cheng2016long,chen2016efficient}.

%%%%% s5
\section{Generalized BDF2, AM2, and IMEX}

Many natural and engineering systems involve energy exchange between different components of the system. Such exchange of energy is often energy conservative, and lead to a skew symmetric term in the model. The rotating fluid model is a well-known example with such a skew symmetric term  \cite{greenspan1968theory}. The Stokes-Darcy system modeling flows in karst aquifers is another example among many others, see for instance  \cite{discacciati2003analysis} and the Appendix.
 
Direct applications of the generalized BDF2 or generalized AM2 method to such problems with a skew-symmetric term lead to solving a non-symmetric matrix problem at each time step. Solving such non-symmetric problems is usually more time consuming than solving their symmetric counterparts. In order to enhance the efficiency, we naturally resort to the implicit-explicit methodology (IMEX). In particular, we propose to treat the skew symmetric term explicitly %, either via Adams-Bashforth method in the AM2 case or via Gear's extrapolation in the BDF2 case, 
so that we only need to solve a symmetric problem if the damping term is symmetric to start with. 
The development and analysis of efficient and stable schemes that combine generalized BDF2/AM2 with IMEX are the focus of this section.

In order to formulate our main results that are applicable to suitable evolutionary PDE models, we introduce a few notations.
Let $H$ be a Hilbert space, and let $V$ be a subspace of $H$. Denoting  $V'$ as the dual space of $V$ induced by the inner product on $H$. We have $V \subset H \subset V'$. 
A linear operator $\mathcal{L}$ from $V$ to $V'$ is called positive on $V$
%$$\mathcal{L}: V \rightarrow V',$$
 if $\exists l_0>0$ s.t. $\langle \mathcal{L}y,y \rangle \ge l_0\|y\|_H^2, \forall y \in V$. 
 A linear operator $\mathcal{L}_s: V \rightarrow V'$ is called {\it skew-symmetric} 
 if $\mathcal{L}_s^*=-\mathcal{L}_s$, where $\mathcal{L}_s^*$ represents the adjoint of $\mathcal{L}_s$. Skew symmetry implies $\langle \mathcal{L}_s y ,y \rangle = 0, \forall y \in V$.
 
We focus on the following type of model with both a dissipative term and a skew-symmetric term.
\begin{equation}\label{skew}
\frac{dy}{dt} + \mathcal{L}y + \mathcal{L}_{s} y = g(t),
\end{equation}
where $f$ is uniformly bounded in time, $\mathcal{L}$ is positive definite, and $\mathcal{L}_s$ is skew symmetric.

Well-known examples that fall into this category include the linear rotating fluid equation \cite{greenspan1968theory} as well as the Stokes-Darcy system for flow in karst aquifer\cite{chen2013efficient,discacciati2009navier}.

 \subsection{Generalized BDF2-IMEX scheme}
 In order to develop efficient and stable algorithms, we follow conventional wisdom,  and treat the ``nice'' dissipative term implicitly for stability, and the skew-symmetric term  explicitly for efficiency. When combined with the generalized BDF2 scheme, we arrive at the following generalized BDF2 based IMEX method for \eqref{skew}:
 \begin{equation}\label{IMEX}
 \frac{3}{2}y_{n+1} - 2y_n + \frac{1}{2}y_{n-1} + h\mathcal{L}(\alpha y_{n+1} + (2-2\alpha)y_{n} + (\alpha-1)y_{n-1}) + h\mathcal{L}_s(2y_n-y_{n-1}) = hg_{n+1},
 \end{equation}
 after utilizing Gear's extrapolation on the skew symmetric term.
Such a hybrid implicit-explicit approach  (IMEX schemes) have been extensively studied\cite{anitescu2004implicit,ascher1995implicit,chen2013efficient,chen2016efficient,frank1997stability,gottlieb2012long,shen2023split,jiang2016optimally,layton2012stability}.

The weak formulation of $\eqref{skew}$, which is required for applications to evolutionary PDEs in general,  can be derived by formally taking the inner product of the equation with a test function $\tilde{y}$ in $V$. This is given as follows: given $g \in H$, seek $y \in L^2(0, T; V), \frac{dy}{dt}\in L^2(0,T; V')$, s.t.
$$\langle \frac{dy}{dt}, \tilde{y} \rangle + \langle \mathcal{L}y,\tilde{y} \rangle + \langle \mathcal{L}_s y, \tilde{y} \rangle = \langle g, \tilde{y} \rangle, \forall \tilde{y} \in V.$$

The weak formulation of the generalized BDF2-IMEX method then takes the form
\begin{equation}\label{IMEX-WeakForm}
    \begin{aligned}
      &  \langle \frac{3}{2}y_{n+1} - 2y_n + \frac{1}{2}y_{n-1} ,\tilde{y} \rangle + h\langle \mathcal{L}(\alpha y_{n+1} + (2-2\alpha)y_{n} + (\alpha-1)y_{n-1}), \tilde{y} \rangle \\
        & + h\langle \mathcal{L}_s(2y_n-y_{n-1}), \tilde{y}   \rangle = h\langle g_{n+1}, \tilde{y} \rangle.
    \end{aligned}
\end{equation}

The main conclusion of this subsection is the following stability result.
\begin{theorem}[Unconditional and uniform-in-time energy bound for the generalized BDF2-IMEX method]
The  generalized BDF2-IMEX method \eqref{IMEX-WeakForm} is unconditionally stable  if the skew symmetric term is dominated by the positive dissipative term, i.e., there exists a positive constant $C_1>0$, s.t.,
 \begin{equation}\label{dominance}    \langle \mathcal{L}y, y\rangle \ge C_1\|\mathcal{L}_s y \|_H^2, \forall y \in V. \end{equation}
Moreover, the $H$ norm of the solutions to this scheme are uniformly bounded for all time if the time-step restriction \eqref{gBDF2-step-restriction} holds.

\end{theorem}
\begin{proof}
For simplicity, we denote
$$ D_{\alpha} y_{n+1} = \alpha y_{n+1} + (2-2\alpha)y_n + (\alpha -1)y_{n-1}.$$

Letting  $\tilde{y} = D_{\alpha} y_{n+1} =\alpha y_{n+1} + (2-2\alpha) y_n + (\alpha - 1)y_{n-1}$ in \eqref{IMEX-WeakForm}, we obtain,
\begin{equation}\label{energy_inequ_with_bdFd_skew_sym}
    \begin{aligned}
    &\|\bV_{n+1}\|_G^2  + \frac{4\alpha-3}{4}\|y_{n+1} -2y_n + y_{n-1}\|_H^2 + h \langle \mathcal{L}(D_{\alpha} y_{n+1}),D_{\alpha} y_{n+1} \rangle  + h\langle \mathcal{L}_s(2y_n-y_{n-1}),D_{\alpha}y_{n+1}\rangle  \\
    =& \|\bV_n\|_G^2 + h \langle g_{n+1},D_{\alpha} y_{n+1}\rangle,
    \end{aligned}
\end{equation}
where $\bV_{n+1} = [y_n,y_{n+1}]^T$,  $\|\bV\|_G^2 = \langle \bV, G \bV \rangle$ with the corresponding G matrix given by
\begin{equation}
G = \frac{1}{4}
    \begin{pmatrix}
   2 \alpha -1 & -2\alpha \\ -2\alpha & 2\alpha +3
    \end{pmatrix},\quad |G| = \frac{4\alpha-3}{16},
\end{equation}
 and we have also utilized the following identity
\begin{equation}
\begin{split}
    & \langle \frac{3}{2}y_{n+1} -2y_n + \frac{1}{2}y_{n-1} , \alpha y_{n+1} + (2-2\alpha)y_n + (\alpha -1)y_{n-1} \rangle   \\
   = &  \|\bV_{n+1}\|_G^2 - \|\bV_n\|_G^2 + \frac{4\alpha-3}{4}\|y_{n+1}-2y_n+y_{n-1}\|_H^2.
\end{split}
\end{equation}
See for instance, Example 6.5 and Theorem 6.6 in \cite{wanner1996solving}.
Notice that the matrix G is positive definite if and only if $\alpha > \frac{3}{4}$.

For the skew symmetric term, it is easy to verify
$$ \langle \mathcal{L}_s(2y_n-y_{n-1}), \alpha y_{n+1} + 2(1-\alpha)y_n - (1-\alpha)y_{n-1}\rangle 
       = \alpha \langle \mathcal{L}_s(\alpha y_{n+1} + 2(1-\alpha)y_n - (1-\alpha)y_{n-1}), y_{n+1}-2y_n + y_{n-1}\rangle.$$
       Hence, 
\begin{equation}
    \begin{split}
        & h\langle \mathcal{L}_s(2y_n-y_{n-1}), \alpha y_{n+1} + 2(1-\alpha)y_n - (1-\alpha)y_{n-1}\rangle \\
       \ge & -\frac{h\alpha^2}{2C_1}\| y_{n+1}-2y_n+y_{n-1}\|_H^2 - \frac{hC_1}{2}\| \mathcal{L}_s(\alpha y_{n+1} + 2(1-\alpha)y_n - (1-\alpha)y_{n-1})\|_H^2, \\
    \end{split}
\end{equation}
Here $ C_1>0$ is the constant in the dominance condition \eqref{dominance}.

For the dissipative term, since $\mathcal{L}$ is positive definite,
% i.e., there exists $l_0>0$, s.t.  $\langle \mathcal{L}v,v \rangle  \ge l_0 \|v\|_H^2$, 
we have,  by Young's inequality,
\begin{equation}\label{dissipative_term_estimate}
\begin{split}
    h \langle \mathcal{L}D_{\alpha}y_{n+1}, D_{\alpha}y_{n+1}\rangle  
    \ge & h l_0\|D_{\alpha} y_{n+1}\|_H^2\\
    = & h l_0\left\{\|(\alpha-1) y_{n+1} + (2-2\alpha)y_n + (\alpha -1)y_{n-1}\|_H^2 + \|y_{n+1}\|_H^2 \right\}\\
    & + 2hl_0\langle (\alpha-1) y_{n+1} + (2-2\alpha)y_n + (\alpha -1)y_{n-1},y_{n+1}\rangle\\
    \ge & \frac{hl_0}{2}\|y_{n+1}\|_H^2 - hl_0(\alpha-1)^2\|y_{n+1} -2 y_n + y_{n-1}\|_H^2.
    \end{split}
\end{equation}

The estimation of the forcing term is straightforward. Indeed, $\forall \varepsilon_2\ge 0$,
\begin{equation}
    \begin{split}
        h\langle g_{n+1},D_{\alpha}y_{n+1} \rangle & =  h \langle g_{n+1},(\alpha-1) y_{n+1} + (2-2\alpha)y_n + (\alpha -1)y_{n-1}\rangle + h\langle g_{n+1},y_{n+1}\rangle \\
        & \le  (\frac{h}{2}+\frac{h}{2\varepsilon_2})\|g\|_{\infty}^2 + \frac{h}{2}(\alpha-1)^2\|  y_{n+1}-2y_n +y_{n-1}\|_H^2 + \frac{h\varepsilon_2}{2}\|y_{n+1}\|_H^2.
    \end{split}
\end{equation}

We split the dissipative term into two equal parts and estimate half of the dissipative term using \eqref{dissipative_term_estimate}, the other half via the dominance condition \eqref{dominance}, and insert the estimate of forcing term into \eqref{energy_inequ_with_bdFd_skew_sym}:

\begin{equation}
    \begin{split}
        & \|\bV_{n+1}\|_G^2  + \frac{4\alpha-3}{4}\|y_{n+1} -2y_n + y_{n-1}\|_H^2 +\frac{hl_0}{4}\|y_{n+1}\|_H^2 \\
    \le& \|\bV_n\|_G^2 + (\frac{h}{2}+\frac{h}{2\varepsilon_2})\|g\|_{\infty}^2 + \left( \frac{hl_0(\alpha-1)^2}{2}+\frac{h(\alpha-1)^2}{2}+\frac{h\alpha^2}{2C_1}\right)\|y_{n+1}-2y_n+y_{n-1}\|_H^2 + \frac{h\varepsilon_2}{2}\|y_{n+1}\|_H^2. \\
    \end{split}
\end{equation}
Setting  $\varepsilon_2 = \frac{l_0}{4}$, we deduce,
\begin{equation}\label{G-estimate-intermediate}
    \begin{split}
        & \|\bV_{n+1}\|_G^2  + \frac{4\alpha-3}{4}\|y_{n+1} -2y_n + y_{n-1}\|_H^2 +\frac{hl_0}{8}\|y_{n+1}\|_H^2 \\
    \le& \|\bV_n\|_G^2 + (\frac{h}{2}+\frac{2h}{l_0})\|g\|_{\infty}^2 + \left( \frac{hl_0(\alpha-1)^2}{2}+\frac{h(\alpha-1)^2}{2}+\frac{h\alpha^2}{2C_1}\right)\|y_{n+1}-2y_n+y_{n-1}\|_H^2. \\
    \end{split}
\end{equation}

Recall the equivalence between the G-norm and standard norm $\|\bV_n\|^2 \triangleq \|y_n\|_H^2 + \|y_{n-1}\|_H^2$, i.e.,  
\begin{equation}
\label{equiv_norm}
C_l \|\bV\|_G^2 \le \|\bV\|^2 \le C_u \|\bV\|_G^2,
\end{equation}
where $C_l = \frac{1}{4} + \frac{\alpha-\sqrt{\alpha^2+1}}{2}, C_u = \frac{1}{4} +  \frac{\alpha+\sqrt{\alpha^2+1}}{2}$.

For the uniform-in-time energy bound, we take $h$ small enough, s.t.,
\begin{equation}\label{gBDF2-step-restriction}
\frac{hl_0(\alpha-1)^2}{2}+\frac{h(\alpha-1)^2}{2}+\frac{h\alpha^2}{2C_1} 
\le \frac{4\alpha-3}{4}.
 \end{equation}
 We then have
\begin{equation}
  \|\bV_{n+1}\|_G^2  +\frac{hl_0}{8}\|y_{n+1}\|_H^2
    \le  \|\bV_n\|_G^2 + (\frac{h}{2}+\frac{2h}{l_0})\|g\|_{\infty}^2.
\end{equation}
Adding $h\delta_2\|y_n\|_H^2$ with $\delta_2=\frac{l_0}{32}$ to both sides of the inequality and utilizing the equivalence of the G-norm and standard norm \eqref{equiv_norm},
%$|\bV_n|^2 \triangleq \|y_n\|_H^2 + \|y_{n-1}\|_H^2$, i.e.,  $C_l \|\bV\|_G^2 \le |\bV|^2 \le C_u \|\bV\|_G^2$
%with $C_l = \frac{1}{4} + \frac{a-\sqrt{a^2+1}}{2}, C_u = \frac{1}{4} +  \frac{a+\sqrt{a^2+1}}{2}$,
we obtain
\begin{equation}\label{G_norm_inequality}
    \|\bV_{n+1}\|_G^2 + \frac{hl_0C_l }{32}\|\bV_{n+1}\|^2_G + \frac{3l_0h}{32}\|y_{n+1}\|_H^2 \le \|\bV_n\|_G^2 + \frac{hl_0}{32}\|y_n\|_H^2+ (\frac{h}{2}+\frac{2h}{l_0})\|g\|_{\infty}^2.
\end{equation}
% correct the constant  
Define
\begin{equation}\label{C_value}
    C_a  = \min (1+\frac{hl_0 C_l}{32},3), E_n = \|\bV_n\|_G^2 + \frac{hl_0}{32}\|y_n\|^2,
 \end{equation}
we deduce
\begin{equation}
    C_aE_{n+1} \le E_n + (\frac{h}{2}+\frac{2h}{l_0})\|g\|_{\infty}^2.
\end{equation}
Iterating the inequality and utilizing the fact that $C_a>1$ leads to 
\begin{equation}
 %   \begin{split}
    E_{n+1} \le   \frac{1}{C_a}E_n + (\frac{h}{2C_a}+\frac{2h}{l_0C_a})\|g\|_{\infty}^2
    \le  \frac{1}{C_a^{n+1}}E_0 +  \frac{1-1/C_a^{n+1}}{C_a-1} (\frac{h}{2}+\frac{2h}{l_0})\|g\|_{\infty}^2.
 %   \end{split}
\end{equation}
This implies a uniform in time energy bound on the solution together with a large-time asymptotic bound independent of the initial data.

For the unconditional stability, we bound the skew symmetry term slightly different as
%$$(Yinqian\ insert)$$
\begin{equation}
    \begin{split}
        h\langle \mathcal{L}_s(2y_n-y_{n-1}), D_{\alpha} y_{n+1} \rangle = & -h \langle 2y_n - y_{n-1}, \mathcal{L}_s D_{\alpha} y_{n+1} \rangle\\
        \ge & -\frac{hC_1}{4}\| \mathcal{L}_s D_{\alpha} y_{n+1}\|^2 - \frac{h}{C_1}\|2y_n-y_{n-1}\|^2\\
        \ge & -\frac{hC_1}{4}\|\mathcal{L}_s D_{\alpha} y_{n+1}\|^2 - \frac{5h}{C_1}(\|y_n\|^2 + \|y_{n-1}\|^2),
    \end{split}
\end{equation}
where we have utilized the skew symmetry of $\mathcal{L}_s$ and Cauchy-Schwarz inequality.  The forcing term can be estimates as %$$(Yinqian\ insert).$$
\begin{equation}
    h\langle g_{n+1}, D_{\alpha}y_{n+1} \rangle \le \frac{hl_0}{4} \|D_{\alpha} y_{n+1}\|^2 + \frac{h}{l_0}\|g_{n+1}\|^2.
\end{equation}

We estimate the dissipative term by splitting it into 1/4, 1/4, 1/2 with 1/4 using the positivity assumption to cancel the $\|D_\alpha y_{n+1}\|^2$ term induced by the forcing term, 1//4 using the dominance assumption to control the $\mathcal{L}_s\|D_\alpha y_{n+1}\|^2$term induced by the skew symmetric term, we have
%simply employ the the equivalence between the G-norm and standard norm $\|\bV_n\|^2 \triangleq \|y_n\|_H^2 + \|y_{n-1}\|_H^2$, i.e.,  $C_l \|\bV\|_G^2 \le \|\bV\|^2 \le C_u \|\bV\|_G^2$
%with $C_l = \frac{1}{4} + \frac{\alpha-\sqrt{\alpha^2+1}}{2}, C_u = \frac{1}{4} +  \frac{\alpha+\sqrt{\alpha^2+1}}{2}$, and we deduce, there exists a positive constant $C$ such that 
\begin{eqnarray}
    \nonumber
 \|\bV_{n+1}\|_G^2  
    &\le& \|\bV_n\|_G^2 + \frac{5h}{C_1}(\|y_n\|^2+\|y_{n-1}\|^2) +\frac{h}{l_0}\|g\|_{\infty}^2 
    \\
    &\le& (1+\frac{5C_u h}{C_1})\|\bV_n\|_G^2 +\frac{h}{l_0}\|g\|_{\infty}^2.
    \end{eqnarray}
The unconditional stability follows in a straightforward manner.

\end{proof}

\noindent{\bf Remark:} One can derive long-time bound as $h\rightarrow 0$ consistent with the long-time bound for \eqref{skew}.

\noindent{\bf Remark:} The time-step restriction \eqref{gBDF2-step-restriction} is mild in the sense that it is independent of spatial discretization.

\subsection{Generalized AM-AB2 scheme}

We now focus on the combination of the generalized AM2 method with IMEX approach to develop stable and efficient numerical algorithms for systems of the form \eqref{skew}.

Following classical Adams-Bashforth approach of treating the skew symmetric term explicitly, we arrive at the following generalized AM-AB2 scheme for \eqref{skew}: 
\begin{equation}\label{gAMB2}
    \begin{aligned}
        &  \langle y_{n+1} - y_{n} ,\tilde{y} \rangle + h\langle \mathcal{L}(\alpha y_{n+1} + (\frac{3}{2}-2\alpha)y_{n} + (\alpha-\frac{1}{2})y_{n-1}), \tilde{y} \rangle \\
        & 
         + h\langle \mathcal{L}_s(\frac{3}{2}y_n-\frac{1}{2}y_{n-1}), \tilde{y}   \rangle = h\langle f_{n+\frac{1}{2}}, \tilde{y} \rangle,
    \end{aligned}
\end{equation}
where $\mathcal{L}$ is a positive definite operator and $\mathcal{L}_s$ is a skew symmetric operator. % s.t. $\mathcal{L}_s^* = -\mathcal{L}_s$. 
The application of this generalized AM-AB2 method to Stokes-Darcy system can be found in \cite{chen2013efficient}. See also \cite{chenwb2016efficient}.

The main conclusion of this subsection is the following stability result.
\begin{theorem}[Unconditional stability and uniform-in-time energy bound for the generalized AM-AB2 method]
    For $\alpha\in (1/2, 1)$, the  generalized AM-AB2 method \eqref{gAMB2} is unconditionally stable under the following conditions:
    \begin{equation}\label{gAMB2-condition}
        \begin{aligned}
            & (a) \text{positive definite leading term: }  \exists l_0 > 0, s.t. \langle \mathcal{L}y,y\rangle \ge l_0 \|y\|_V^2, \forall y \in V \subset H \subset V',\\
            & (b) \text{skew symmetric: } \mathcal{L}_s^*=-\mathcal{L}_s,,\\
            & (c) \text{weak skew symmetric term:}\  \exists  \theta \in (0,1), s.t.,  |\langle \mathcal{L}_sy,z\rangle| \lesssim \|z\|_V^{\frac{1}{2}} \|y\|_H^{\theta}\|y\|_V^{1-\theta},\forall y, z\in V.
        \end{aligned}
    \end{equation}
Moreover, the $H$ norm of the solutions to scheme \eqref{gAMB2} are uniformly bounded in time if the time-step restriction \eqref{gAMB2-step-restriction} holds.

% the following step size restriction holds:
%\begin{equation}
 %      h\le \frac{(\alpha-\alpha_1-\alpha_2)C_1}{12},
%\end{equation}
%where $\alpha$ is the parameter in \eqref{GAM2}, s.t. $\frac{1}{2}<\alpha<1$, and $\alpha_1 \triangleq |\frac{3}{2}-2\alpha|, \alpha_2 \triangleq |\alpha-\frac{1}{2}|.$
\end{theorem}
Note that $A \lesssim B$ means that $A$ is less than or equal to a contant multiply of $B$.
\begin{proof}

Setting the test function $\tilde{y} = y_{n+1}$ in \eqref{gAMB2}, we have
\begin{equation}\label{GAM2_eq1}
\begin{split}
 &     \|y_{n+1}\|_H^2 + \|y_{n+1}-y_n\|_H^2 + 2h\langle \mathcal{L}(D_{\alpha}y_{n+1}) , y_{n+1}\rangle  \\
    & =  \|y_n\|_H^2 + 2h\langle f_{n+\frac{1}{2}}, y_{n+1} \rangle - 2h\langle \mathcal{L}_s(\frac{3}{2}y_n-\frac{1}{2}y_{n-1}),y_{n+1}\rangle,
\end{split}
\end{equation}
where $D_{\alpha}y_{n+1} = \alpha y_{n+1} + (\frac{3}{2}-2\alpha)y_n + (\alpha-\frac{1}{2})y_{n-1}.$

Thanks to the positivity of $\mathcal{L}$ and the Cauchy-Schwarz inequality, we have
\begin{equation}\label{GAM2_diss}
        2h\langle \mathcal{L}(D_{\alpha}y_{n+1}) , y_{n+1}\rangle \ge (2\alpha-(\alpha_1 + \alpha_2))h\langle \mathcal{L}y_{n+1},y_{n+1}\rangle - \alpha_1h \langle \mathcal{L}y_n,y_n\rangle - \alpha_2h \langle \mathcal{L}y_{n-1},y_{n-1}\rangle,
\end{equation}
where 
\begin{equation}
   \alpha_1 = |\frac{3}{2}-2\alpha|, \alpha_2 = |\alpha-\frac{1}{2}|.
\end{equation}
Let
\begin{equation}
         \beta_3 = \alpha_1 + \alpha_2, \beta_1 = 2\alpha-\beta_3, \beta_2 = \frac{\beta_1 +\beta_3}{2}.
\end{equation}
Note that $\alpha \in (\frac{1}{2},1)$ implies  $\beta_1>\beta_2>\beta_3>0,\alpha-\beta_3=\alpha-(\alpha_1+\alpha_2)>0.$

For the skew symmetric term,condition  $(c)$ implies
\begin{equation}\label{GAM2_skew}
\begin{split}
    & - 2h\langle \mathcal{L}_s(\frac{3}{2}y_n-\frac{1}{2}y_{n-1}),y_{n+1}\rangle\\ 
    \le & \frac{\varepsilon h}{2}\|y_{n+1}\|_V^2 + \varepsilon h(\|y_n\|_V^2 + \|y_{n-1}\|_V^2) + Ch(\|y_n\|_H^2 + \|y_{n-1}\|_H^2).
\end{split}
\end{equation}

The forcing term can be bounded as follows
$$2h\langle f_{n+\frac{1}{2}}, y_{n+1} \rangle \le \frac{\varepsilon h}{2}\|y_{n+1}\|_V^2 + Ch\|f_{n+\frac{1}{2}}\|_{V'}^2.$$

Hence we have
\begin{equation}\label{GAM2_eq2}
    \begin{split}
        & \|y_{n+1}\|_H^2 + \|y_{n+1}-y_n\|_H^2-\|y_n\|_H^2 + (2\alpha-(\alpha_1+\alpha_2))h\langle \mathcal{L}y_{n+1},y_{n+1}\rangle \\
        \le & \alpha_1 h\langle \mathcal{L}y_n,y_n\rangle + \alpha_2h\langle \mathcal{L}y_{n-1},y_{n-1}\rangle + Ch (\|y_n\|_H^2 + \|y_{n-1}\|_H^2)\\
        & + \varepsilon h \|y_{n+1}\|_V^2 + 
        \varepsilon h(\|y_n\|_V^2 + \|y_{n-1}\|_V^2) + 
        Ch\|f_{n+\frac{1}{2}}\|_{V'}^2 .\\
    \end{split}
\end{equation}

Notice that 
$$l_0\|y\|_V^2 \le \langle \mathcal{L}y,y\rangle \le l_u\|y\|_V^2.$$

By choosing $\varepsilon$ small enough, s.t.
$$2\alpha-(\alpha_1+\alpha_2) - \frac{\varepsilon}{l_0}>0,$$
we deduce, possible with different constants $C$,
\begin{equation}
    \begin{split}
      &  \|y_{n+1}\|_H^2 + (2\alpha -(\alpha_1 + \alpha_2)-\frac{\varepsilon}{l_0})h\langle \mathcal{L}y_{n+1},y_{n+1}\rangle \\
      \le & \|y_n\|_H^2 + (\alpha_1+\frac{\varepsilon}{l_0})h\langle \mathcal{L}y_n,y_n\rangle + (\alpha_2+\frac{\varepsilon}{l_0})h\langle \mathcal{L}y_{n-1},y_{n-1}\rangle + Ch\|f_{n+\frac{1}{2}}\|_{V'}^2 + Ch(\|y_n\|_H^2 + \|y_{n-1}\|_H^2).
    \end{split}
\end{equation}

Adding $\gamma h\langle \mathcal{L}y_n,y_n\rangle $ to both sides to obtain
\begin{equation}
    \begin{split}
         &  \|y_{n+1}\|_H^2 + (2\alpha -(\alpha_1 + \alpha_2)-\frac{\varepsilon}{l_0})h\langle \mathcal{L}y_{n+1},y_{n+1}\rangle + \gamma h\langle \mathcal{L}y_n,y_n\rangle \\
      \le & \|y_n\|_H^2 + (\alpha_1+\frac{\varepsilon}{l_0}+\gamma)h\langle \mathcal{L}y_n,y_n\rangle + (\alpha_2+\frac{\varepsilon}{l_0} )h\langle \mathcal{L}y_{n-1},y_{n-1}\rangle \\
      & + Ch\|f_{n+\frac{1}{2}}\|_{V'}^2 + Ch(\|y_n\|_H^2 + \|y_{n-1}\|_H^2).
    \end{split}
\end{equation}

Letting $\gamma = \alpha_2 + \frac{\varepsilon}{l_0}$ and adding $h\|y_n\|_H^2$ to both sides, we have
\begin{equation}
    \begin{split}
         &  \|y_{n+1}\|_H^2 + h\|y_n\|_H^2 + (2\alpha -(\alpha_1 + \alpha_2)-\frac{\varepsilon}{l_0})h\langle \mathcal{L}y_{n+1},y_{n+1}\rangle + (\alpha_2 + \frac{\varepsilon}{l_0}) h\langle \mathcal{L}y_n,y_n\rangle \\
      \le & (1+Ch) \|y_n\|_H^2 + Ch\|y_{n-1}\|_H^2 \\
      & +  (\alpha_1+\alpha_2+\frac{2\varepsilon}{l_0})h\langle \mathcal{L}y_n,y_n\rangle + (\alpha_2+\frac{\varepsilon}{l_0} )h\langle \mathcal{L}y_{n-1},y_{n-1}\rangle 
       + Ch\|f_{n+\frac{1}{2}}\|_{V'}^2 .
    \end{split}
\end{equation}

Choose $\varepsilon$ small enough so that
$$2\alpha - (\alpha_1 + \alpha_2)-\frac{\varepsilon}{l_0} \ge \alpha_1 + \alpha_2 + \frac{2\varepsilon}{l_0}.$$
This is doable since  $\alpha > \alpha_1+\alpha_2$, which is a consequence of $\alpha \in (\frac{1}{2},1)$.

Denoting $$E_{n+1} = \|y_{n+1}\|_H^2 + h\|y_n\|_H^2 + (2\alpha -(\alpha_1 + \alpha_2)-\frac{\varepsilon}{l_0})h\langle \mathcal{L}y_{n+1},y_{n+1}\rangle + (\alpha_2 + \frac{\varepsilon}{l_0}) h\langle \mathcal{L}y_n,y_n\rangle ,$$
we have
$$E_{n+1} \le (1+Ch)E_n  + Ch \|f_{n+\frac{1}{2}}\|_{V'}^2,$$
and the unconditional stability follows in a straightforward manner.

For the purpose of uniform-in-time energy bound of the solutions to the generalized AM-AB2 scheme, we revisit the estimate of the skew-symmetric term. Utilizing the skew symmetry, we deduce
$$\langle \mathcal{L}_s (\frac{3}{2}y_n-\frac{1}{2}y_{n-1}),y_{n+1}\rangle = -\langle \mathcal{L}_s(y_{n+1}-y_n),y_{n+1}\rangle + \frac{1}{2}\langle \mathcal{L}_s(y_n-y_{n-1}),y_{n+1}\rangle.$$

Hence, for $\forall \varepsilon_1,\varepsilon_2>0$, after taking advantage of the H{\"o}lder inequality and the dominance of the skew-symmetric term by the dissipative term, i.e. $(c)$, and the equivalence between the $V$ norm and $\langle \mathcal{L}y,y\rangle^{\frac{1}{2}}$, i.e. $(a)$, we have
\begin{equation}
    \begin{aligned}
       & 2|\langle \mathcal{L}_s(\frac{3}{2}y_n-\frac{1}{2}y_{n-1}),y_{n+1}\rangle| \\
        \lesssim & \langle \mathcal{L}y_{n+1},y_{n+1}\rangle^{\frac{1}{2}}(\langle \mathcal{L}(y_{n+1}-y_n),y_{n+1}-y_n\rangle^{\frac{1}{4}}\|y_{n+1}-y_n\|_H^{\frac{1}{2}} +\langle \mathcal{L}(y_{n}-y_{n-1}),y_{n}-y_{n-1}\rangle^{\frac{1}{4}}\|y_{n}-y_{n-1}\|_H^{\frac{1}{2}} )\\
        \le & \varepsilon_1\langle \mathcal{L}y_{n+1},y_{n+1}\rangle + \varepsilon_2 \langle \mathcal{L}y_{n},y_{n}\rangle + \varepsilon_2 \langle \mathcal{L}y_{n-1},y_{n-1}\rangle + C_s(\|y_{n+1}-y_n\|_H^2 + \|y_n-y_{n-1}\|_H^2).
    \end{aligned}
\end{equation}
Here ``$\lesssim$'' means ``$\le$'' modulo a prefactor. 

Combining this estimate with the dissipative term estimate and the forcing term estimate, we deduce
\begin{equation}
    \begin{aligned}
     &   \|y_{n+1}\|_H^2 + \|y_{n+1}-y_n\|_H^2 + (2\alpha-(\alpha_1+\alpha_2)-\varepsilon_1)h\langle \mathcal{L}y_{n+1},y_{n+1}\rangle\\
     \le & \|y_n\|_H^2 + (\alpha_1+\varepsilon_2)h\langle \mathcal{L}y_n,y_n\rangle + (\alpha_2+\varepsilon_2)h\langle \mathcal{L}y_{n-1},y_{n-1}\rangle \\
     & + C_sh(\|y_{n+1}-y_n\|_H^2 + \|y_n-y_{n-1}\|_H^2) + \varepsilon_3h\langle \mathcal{L}y_{n+1},y_{n+1}\rangle + C_fh\|f_{n+\frac{1}{2}}\|_{V'}^2,
    \end{aligned}
\end{equation}
where we have utilized a slightly different estimate on the forcing term
\begin{equation}
    \begin{aligned}
        2|\langle f_{n+\frac{1}{2}},y_{n+1}\rangle| & \lesssim \langle \mathcal{L}y_{n+1},y_{n+1}\rangle^{\frac{1}{2}}\|f_{n+\frac{1}{2}}\|_{V'}\\
        & \le \varepsilon_3\langle \mathcal{L}y_{n+1},y_{n+1}\rangle + C_f\|f_{n+\frac{1}{2}}\|_{V'}^2.
    \end{aligned}
\end{equation}

Let $\delta,\varepsilon_1,\varepsilon_2,\varepsilon_3$ be small enough and positive, s.t.
\begin{equation}
2\alpha-(\alpha_1+\alpha_2)-\varepsilon_1-\varepsilon_3-\delta > (\alpha_1+\varepsilon_2) + (\alpha_2 + 2\varepsilon_2).
\label{small_parameter_condition}
\end{equation}
Notice, the above inequality is essentially the same as $\alpha-(\alpha_1+\alpha_2)>0$, which is equivalent to $\alpha\in (1/2, 1)$.

Hence it is always possible to choose $\delta,\varepsilon_1,\varepsilon_2,\varepsilon_3>0$. Therefore, we have, after utilizing the equivalent norm on $V$ and the fact that the $V$ norm dominates the $H$ norm, i.e.
$$C_i\|y\|_H^2 \le \|y\|_V^2, \forall y \in V \text{ (imbedding)},$$

\begin{equation}
    \begin{aligned}
     &   (1+C_i\delta l_0h)\|y_{n+1}\|_H^2 + \|y_{n+1}-y_n\|_H^2 \\
     &+ (2\alpha-(\alpha_1 +\alpha_2)-\varepsilon_1-\varepsilon_3-\delta)h\langle \mathcal{L}y_{n+1},y_{n+1}\rangle + (\alpha_2 + 2\varepsilon_2)\langle \mathcal{L}y_n,y_n\rangle\\
     \le & \|y_n\|_H^2 + (\alpha_1+\alpha_2+3\varepsilon_2)h\langle \mathcal{L}y_n,y_n\rangle + (\alpha_2+\varepsilon_2)h\langle \mathcal{L}y_{n-1},y_{n-1}\rangle \\
     & + C_sh(\|y_{n+1}-y_n\|_H^2 + \|y_n-y_{n-1}\|_H^2) + C_fh\|f_{n+\frac{1}{2}}\|_{V'}^2.
    \end{aligned}
\end{equation}

Choose the time-step $h$ small enough so that
\begin{equation}\label{gAMB2-step-restriction}
\begin{aligned}
   & h<\frac{1}{3C_s}, \quad  \frac{2\alpha-(\alpha_1+\alpha_2)-\varepsilon_1-\varepsilon_3-\delta}{1+C_i\delta l_0 h}>(\alpha_1+\alpha_2)+3\varepsilon_2,\\
  &\frac{\alpha_2+2\varepsilon_2}{1+C_i\delta l_0 h}>\alpha_2 + \varepsilon_2,\quad \frac{1-C_sh}{1+C_i\delta l_0 h}>C_sh,
\end{aligned}
\end{equation}

we have
\begin{equation}
    \begin{aligned}
    &    (1+C_i\delta l_0 h)\left\{
    \|y_{n+1}\|_H^2 + \frac{1-C_sh}{1+C_i\delta l_0 h}\|y_{n+1}-y_n\|_H^2 \right.\\
    & \left. +  \frac{2\alpha -(\alpha_1+\alpha_2)-\varepsilon_1-\varepsilon_3-\delta}{1+C_i\delta l_0 h}h\langle \mathcal{L}y_{n+1},y_{n+1}\rangle + \frac{\alpha_2+2\varepsilon_2}{1+C_i\delta l_0 h}h\langle \mathcal{L}y_n,y_n\rangle    \right\}\\
    \le & \|y_n\|_H^2 + (\alpha_1+\alpha_2+3\varepsilon_2)h\langle \mathcal{L}y_n,y_n\rangle + (\alpha_2+\varepsilon_2)h\langle \mathcal{L}y_{n-1},y_{n-1}\rangle \\
    & + C_sh\|y_{n}-y_{n-1}\|_H^2 + C_fh\|f_{n+\frac{1}{2}}\|_{V'}^2.
    \end{aligned}
\end{equation}

Let 
\begin{equation}
    \begin{aligned}
        E_n = & \|y_n\|_H^2 + \frac{1-C_sh}{1+C_i\delta l_0 h}\|y_n-y_{n-1}\|_H^2 \\
        & + \frac{2\alpha-(\alpha_1+\alpha_2)-\varepsilon_1-\varepsilon_3-\delta}{1+C_i \delta l_0 h}h\langle \mathcal{L}y_n,y_n\rangle + \frac{\alpha_2 + 2\varepsilon_2}{1+C_i\delta l_0 h}h\langle \mathcal{L}y_{n-1},y_{n-1}\rangle.
    \end{aligned}
\end{equation}

We have, by our choice of $\varepsilon_1,\varepsilon_2,\varepsilon_3,\delta$ and $h$,
$$(1+C_i\delta l_0 h)E_{n+1} \le E_n + C_f h \|f_{n+\frac{1}{2}}\|_{V'}^2.$$
The uniform-in-time energy bound follows.

\end{proof}

%%%%%
\section{Numerical experiments}

In this section we present a few simple numerical tests verifying the second order accuracy in short and long time of the generalized BDF2 and AM2 schemes.
 We also demonstrate that the time-step restrictions are needed in the case of IMEX schemes for problems with a skew-symmetric term.
In addition, we point out that the generalized scheme could be more or less accurate, or more or less stable than the classical schemes, depending on the parameter values. 

\subsection{Numerical results for the accuracy test}
\subsubsection{Accuracy test for the generalized BDF2 method}
Consider the following scalar damped-driven ODE with a non-periodic forcing term 
\begin{equation}
\label{damped_driven}
    y' + 10y = sin(t) + cos(\sqrt{2}t), \quad y_0 = 1,
\end{equation}
and apply the generalized BDF2 method \eqref{generalized_BDF2_general_form} with different values of $\alpha$ to the end time $T_{end} = 1 $, to obtain its relative error and convergence rate. Note that in the numerical part, $y_1$ is obtained by substituting the exact solution $y_1 = y(h)$ to start the multi-step method. We could also use a one-step method to initiate the two-step scheme.

%% To answer the review 'too tight to get the conclusion, more values of $\alpha$ are put here.
\begin{figure}[htbp]
    \centering  
    \includegraphics[width=0.6\textwidth]{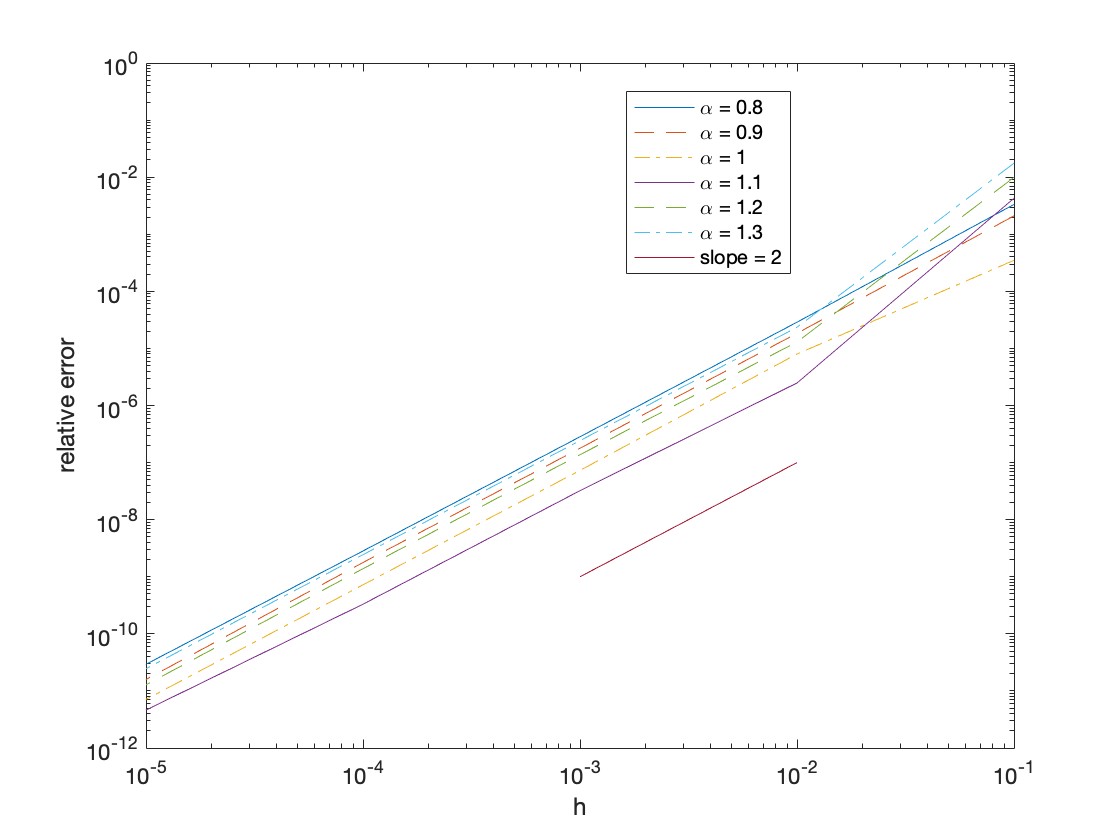}
    \caption{Loglog plot for the relative error of the generalized BDF2 method}
\end{figure}

\begin{table}[htbp]
  \caption{Relative error and order of accuracy for the generalized BDF2 method}
    \centering
    \begin{tabular}{lllllll}
    \noalign{\smallskip}
\hline
\noalign{\smallskip}
{}&\multicolumn{2}{c}{$\alpha = 0.8$} &  \multicolumn{2}{c}{$\alpha = 0.9$}& \multicolumn{2}{c}{$\alpha = 1$}\\
\noalign{\smallskip}
\hline
\noalign{\smallskip}
Time step  &Relative error &  Order & Relative error & Order  &Relative error &  Order\\
\noalign{\smallskip}
\hline
\noalign{\smallskip}
1e-01  &   3.3324e-03&  & 2.1002e-03 & &  3.4969e-04 & \\

1e-02  &  2.8796e-05&  2.0634 & 1.8399e-05&2.0575 &  7.9740e-06  &   1.6420\\

1e-03   & 2.8348e-07 &2.0068 &  1.7825e-07 &   2.0138&  7.3022e-08 & 2.0382\\

1e-04  &2.8301e-09  & 2.0007 &  1.7764e-09 &    2.0015 &   7.2271e-10 & 2.0045\\

1e-05  & 2.9344e-11 &  1.9843 & 1.5853e-11 & 2.0494&   7.1039e-12 & 2.0075\\
\noalign{\smallskip}
\hline
\end{tabular}
\end{table}

\begin{table}[htbp]
  \caption{Relative error and order of accuracy for the generalized BDF2 method}
    \centering
    \begin{tabular}{lllllll}
    \noalign{\smallskip}
\hline
\noalign{\smallskip}
{}&\multicolumn{2}{c}{$\alpha = 1.1$} &  \multicolumn{2}{c}{$\alpha = 1.2$}& \multicolumn{2}{c}{$\alpha = 1.3$}\\
\noalign{\smallskip}
\hline
\noalign{\smallskip}
Time step  &Relative error &  Order & Relative error & Order  &Relative error &  Order\\
\noalign{\smallskip}
\hline
\noalign{\smallskip}
1e-01  &  4.2729e-03 &  &  9.9380e-03& &  1.7597e-02  & \\

1e-02  & 2.4787e-06 &  3.2364 &  1.2959e-05  &  2.8847&  2.3468e-05 &  2.8750\\

1e-03   & 3.2209e-08 &1.8863 & 1.3744e-07 &  1.9745 &  2.4268e-07  & 1.9854\\

1e-04  &3.3073e-10  & 1.9885  &   1.3844e-09  &    1.9969 &   2.4381e-09 &1.9980\\

1e-05  &  4.6301e-12 & 1.8539 &  1.2801e-11 &  2.0340&  2.4527e-11
 & 1.9974\\
\noalign{\smallskip}
\hline
\end{tabular}
\end{table}

Tables 1 and 2 clearly demonstrate the second order accuracy of the generalized BDF2 scheme.
Figure 1 indicates that the generalized BDF2 method may be more accurate (smaller relative error) than the classical BDF2 method ($\alpha = 1$) when $\alpha$ is slightly bigger than the classical value. However, the error becomes larger when $\alpha=1.3$. %{\color{red}We speculate that this gain of accuracy can be lost with huge $\alpha$. Please numerically verify this with $\alpha=1000$}

\subsubsection{Accuracy test for the generalized AM2 method}
Here we apply the generalized AM2 method \eqref{generalized_AM2_general_form} to the same damped driven model \eqref{damped_driven} in order to test the accuracy of the scheme.

% new
\begin{figure}[ht]
    \centering   \includegraphics[width=0.6\textwidth]{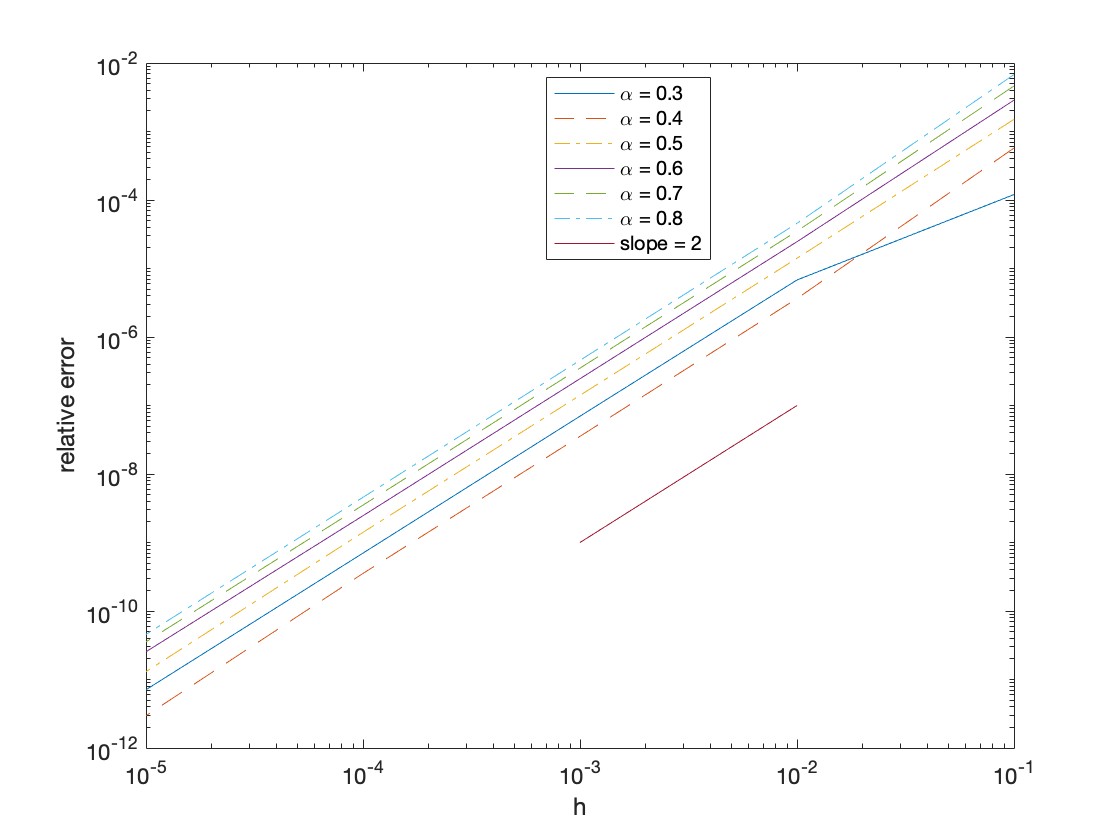}
    \caption{Loglog plot for the relative error of the generalized AM2 method}
\end{figure}

\begin{table}[htbp]
\caption{Relative error and order of accuracy for the generalized AM2 method}
    \centering
   \begin{tabular}{lllllll}
   \noalign{\smallskip}
\hline
\noalign{\smallskip}
{}&\multicolumn{2}{c}{$\alpha = 0.3$} &  \multicolumn{2}{c}{$\alpha = 0.4$}&\multicolumn{2}{c}{$\alpha = 0.5$}\\
\noalign{\smallskip}
\hline
\noalign{\smallskip}
Time step  &Relative error &  Order & Relative error & Order &Relative error &  Order  \\
\noalign{\smallskip}
\hline
\noalign{\smallskip}

1e-01  &    1.2018e-04 &  &5.7368e-04 &  &  1.5054e-03& \\

1e-02  & 6.8001e-06&  1.2473 &3.6451e-06 & 2.1970&1.4116e-05 &2.0279\\

1e-03   &  6.9839e-08&  1.9884& 3.5473e-08 & 2.0118& 1.4079e-07&2.0012\\

1e-04  &  6.9993e-10 &   1.9990& 3.5382e-10 &  2.0011&1.4077e-09&2.0001\\

1e-05  & 7.0959e-12 & 1.9941& 2.9631e-12 & 2.0770&1.3019e-11&2.0339\\
\noalign{\smallskip}
\hline
\end{tabular}
\end{table}

\begin{table}[htbp]
\caption{Relative error and order of accuracy for the generalized AM2 method}
    \centering
   \begin{tabular}{lllllll}
   \noalign{\smallskip}
\hline
\noalign{\smallskip}
{}&\multicolumn{2}{c}{$\alpha = 0.6$} &  \multicolumn{2}{c}{$\alpha = 0.7$}&\multicolumn{2}{c}{$\alpha = 0.8$}\\
\noalign{\smallskip}
\hline
\noalign{\smallskip}
Time step  &Relative error &  Order & Relative error & Order &Relative error &  Order  \\
\noalign{\smallskip}
\hline
\noalign{\smallskip}

1e-01  &   2.8702e-03  &  &4.6227e-03&  & 6.7831e-03 &\\

1e-02   &2.4613e-05 & 2.0667&3.5137e-05&2.1191 & 4.5686e-05& 2.1716\\

1e-03   &   2.4610e-07 &  2.0001&3.5142e-07&2.0000 & 4.5675e-07 & 2.0001\\

1e-04  &   2.4614e-09 &   1.9999&3.5152e-09&1.9999 & 4.5690e-09 & 1.9999\\

1e-05  &  2.5419e-11 & 1.9860&3.5472e-11&1.9961& 4.5540e-11 & 2.0014\\
\noalign{\smallskip}
\hline
\end{tabular}
\end{table}

Tables 3 and 4 clearly demonstrates the second order accuracy of the generalized AM2 scheme.
Figure 2 suggests that generalized AM2 can be more accurate than the classical AM2 method when $\alpha$ is slighly below the threshold value of 0.5 at $T_{end}=1$.
Figure 3 suggests that the generalized AM2 method could be more accurate than the classical AM2 method, even at very large time  ($T_{end}=10^4$) for $\alpha$ slightly below the threshold value. Note that the scheme is no longer A-stable for these values of the parameter $\alpha$.

\begin{figure}[ht]
    \centering   \includegraphics[width=0.7\textwidth]{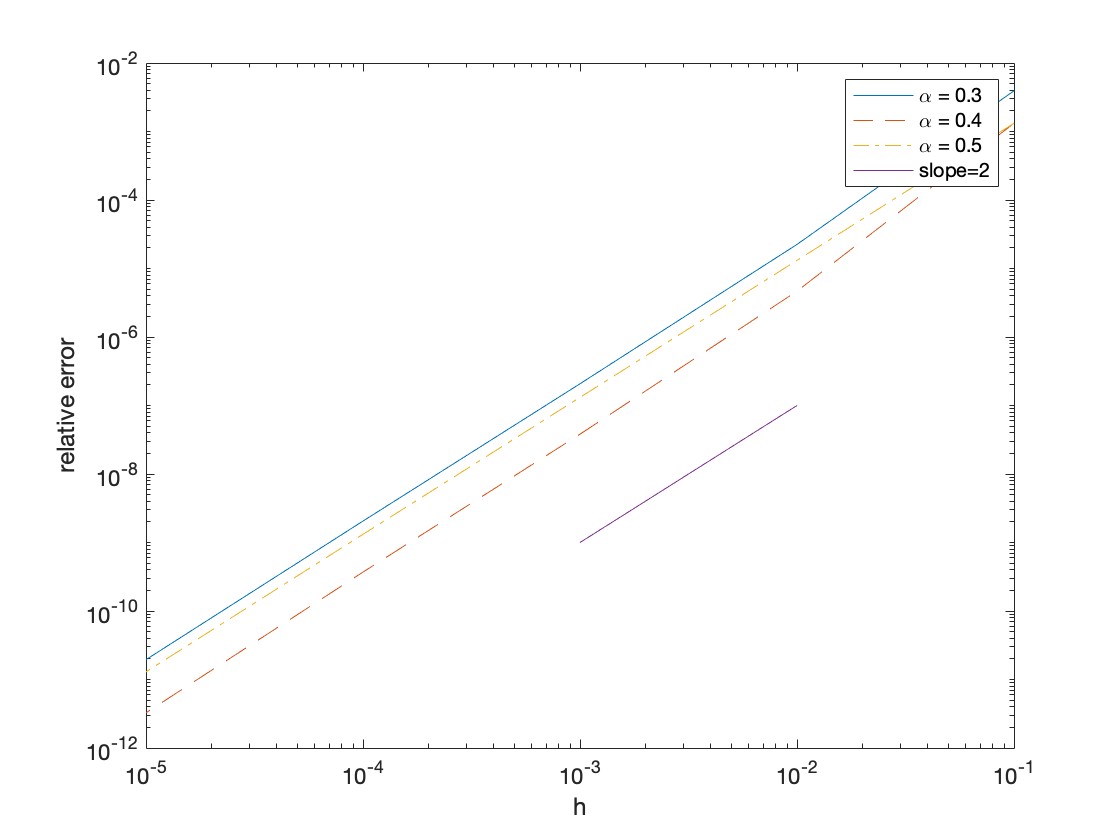}
    \caption{Loglog plot for the relative error of the generalized AM2 method with $T_{end}=10^4$}
\end{figure}

%Note that the theoretical analysis about long-time energy stability for the Stokes-Darcy system can be found in \cite{chen2013efficient,chenwb2016efficient}.

%6.2
\subsection{Numerical results for the energy stability test}

\subsubsection{Numerical experiment on the stabilizing effect}

%Recall that large $\alpha$ corresponds to large stability domain for the generalized BDF2 scheme.
We first show that the time-step restriction is needed for the IMEX BDF2 scheme \eqref{IMEX} when applied to a damped-driven model with a skew symmetric term.
We then demonstrate the uniform-in-time second order accuracy.

For simplicity, we utilize the following  two-dimensional damped-driven ODE system with a non-periodic forcing term 

\begin{equation}
\label{damped_driven_skew}
    \by' + 10\by + \begin{pmatrix}
        -y_{i} \\ y_{r}
    \end{pmatrix}= \begin{pmatrix}
        sin(t)+cos(\sqrt{2}t)\\ 0
    \end{pmatrix},
\end{equation}
with notation $\by = \begin{pmatrix}
    y_{r}\\y_{i}
\end{pmatrix}$ and initial condition as $\by_0 =\begin{pmatrix}
    1\\0
\end{pmatrix} $.
Notice that the third term on the left hand side is skew symmetric.
We conduct numerical experiments with time step size $h=0.1,10$ and the end time $T_{end} = 100,1\times 10^4$.

The {first component of the }approximation solution obtained from the generalized IMEX BDF2 method with different parameter $\alpha$ and 
the {first component of the} exact solution are plotted in Figure 4.

\begin{figure}[ht]
\centering
%\subfigure[$h=0.1, T=100$]{
  \begin{minipage}[t]{0.49\linewidth} 
  \centering
    \includegraphics[width=1\textwidth]{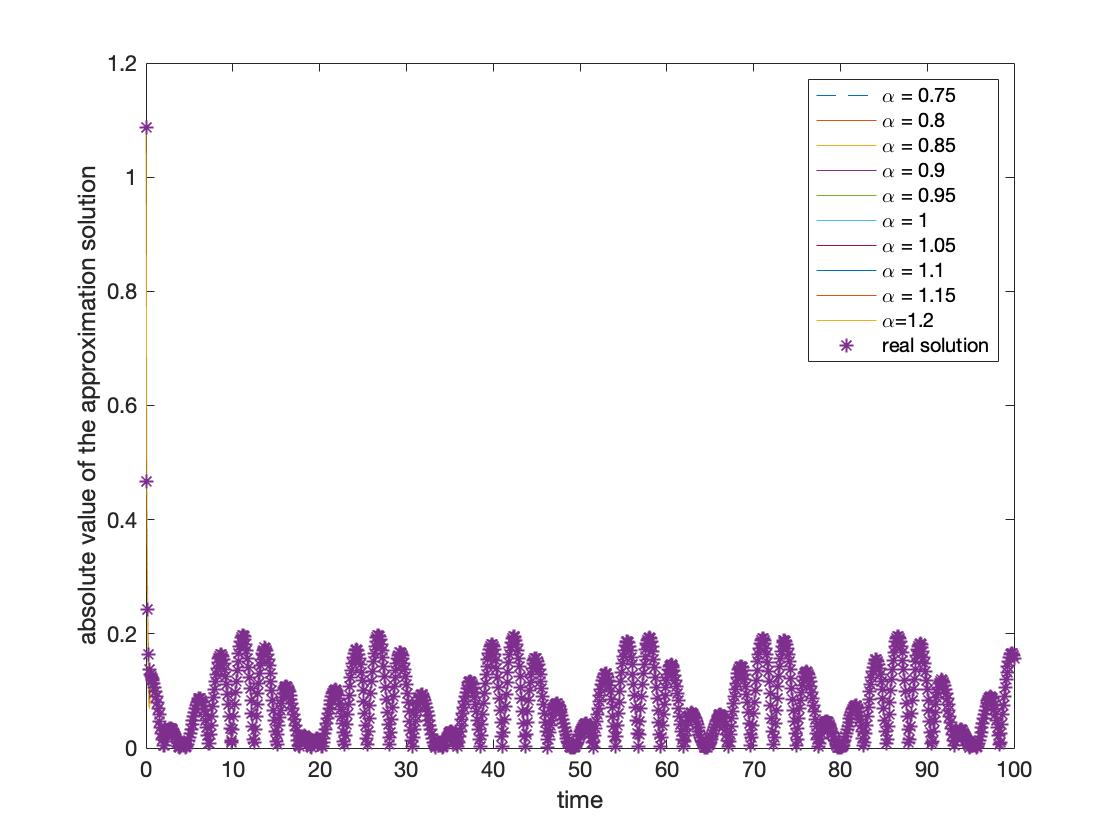}
   % \captionsetup{font={footnotesize},justification=raggedright}
  %  \caption{$h=0.1, T=100$}
    \end{minipage}%
%}
%\subfigure[$h=10, T=10000$]{
\begin{minipage}[t]{0.49\linewidth}
    \includegraphics[width=1\textwidth]{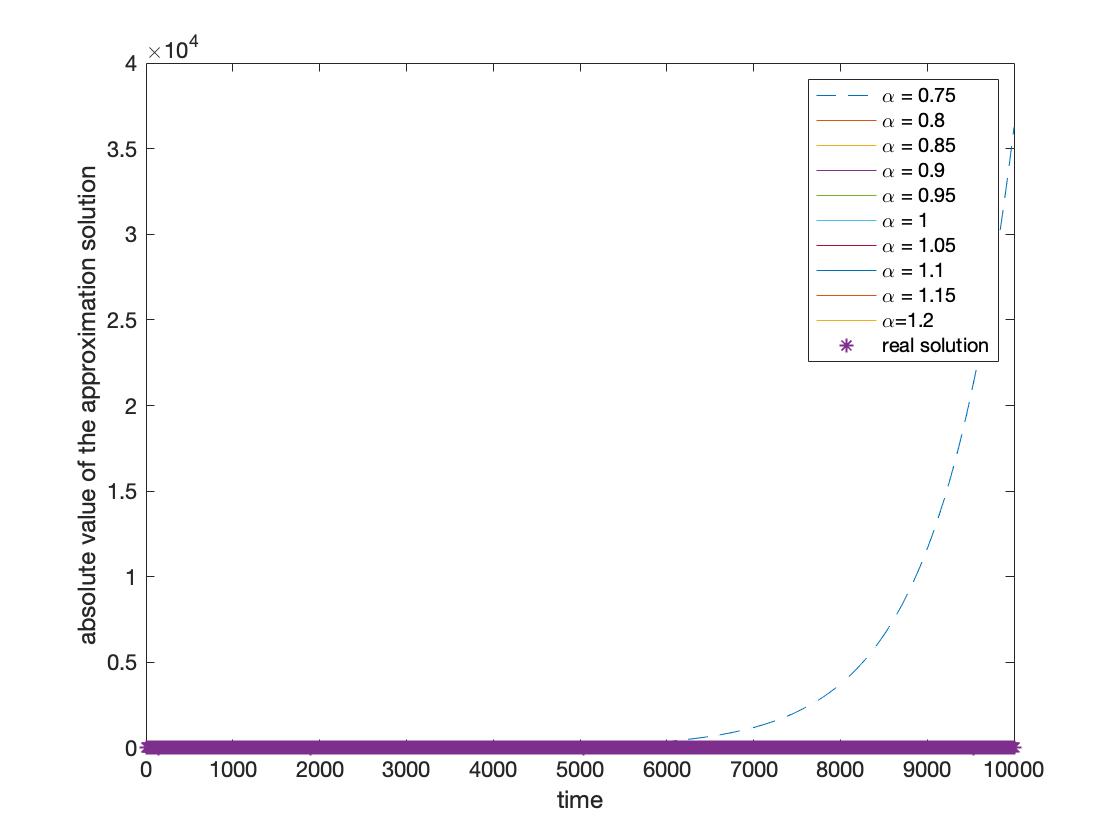}
   %  \captionsetup{font={footnotesize}}
 %   \caption{$h=10, T=10000$}
    \end{minipage}%
%}
\caption{Approximation solution of the generalized IMEX BDF2 method with different parameter $\alpha$. Left panel $h=0.1, T=100$. Right panel $h=10, T=10000$. }
\end{figure}

We observe that the approximation solution is good for small time-step size ($h=0.1$) in the sense that it is almost indistinguishable from the exact solution in the left panel in  figure 4.  However, the right panel in figure 4 indicates that the numerical solutions may blow up (grow exponentially) if the mild time-step size restriction does not hold, in the case of $\alpha = \frac{3}{4}$. This suggests that the time-step restriction imposed in our theory is in fact needed.

When numerical solutions tend to grow exponentially, Figures 5 indicate that larger values of the parameter $\alpha$ may lead to enhanced stability in the sense that the rate of growth may be smaller, when applied to the following damped-driven 2D ODE system with a skew-symmetric term

\begin{equation}
\label{damped_driven_skew2}
 \by' + 10\by + 4\begin{pmatrix}
        -y_{i} \\ y_{r}
    \end{pmatrix}= \begin{pmatrix}
        sin(t)+cos(\sqrt{2}t)\\ 0
    \end{pmatrix}, \quad \by = \begin{pmatrix}
        y_r \\ y_i
    \end{pmatrix}, \by_0 = \begin{pmatrix}
        1 \\ 0
    \end{pmatrix},
\end{equation}
 with $h=100,200,$ and $ T = 2 \times 10^4, 1\times 10^5$.

\begin{figure}[ht]
\centering
%\subfigure[$h=100, T=2\times 10^4$]{
 \begin{minipage}[t]{0.49\linewidth}
    \centering   \includegraphics[width=1\textwidth]{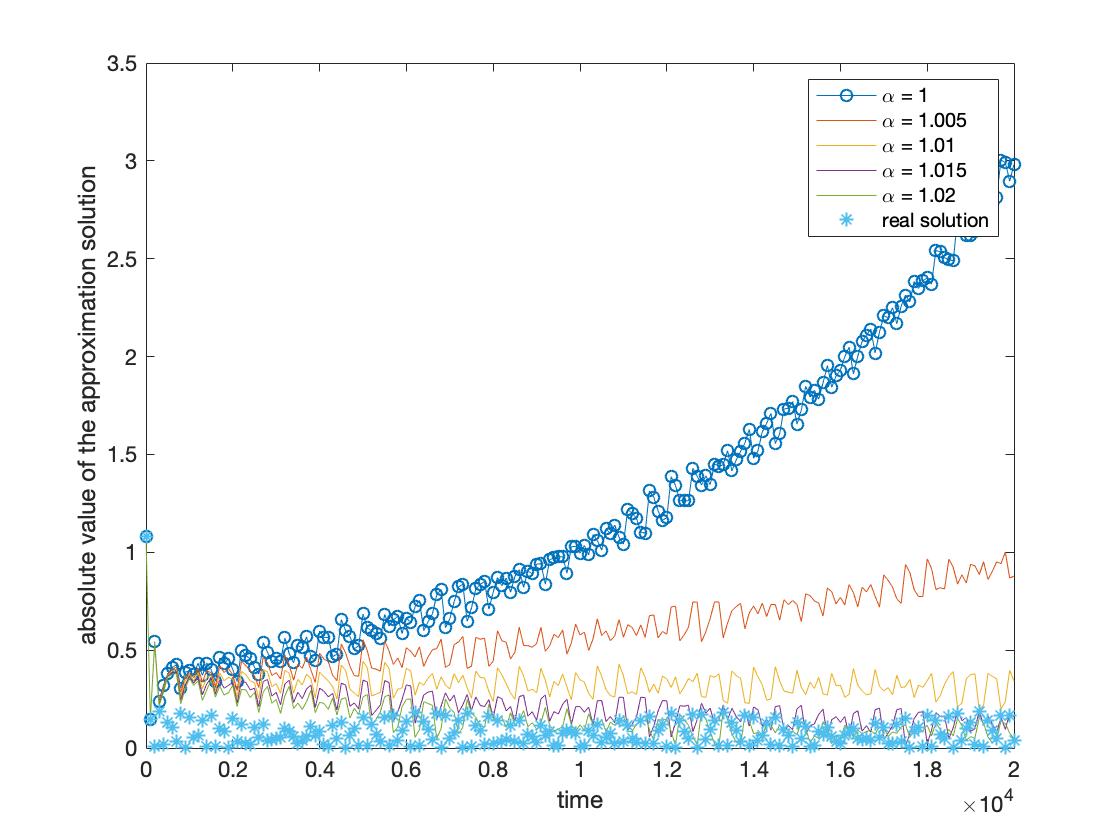}
    % \captionsetup{font={footnotesize}}
 %   \caption{$h=100, T=2\times 10^4$}
    \end{minipage}%
%}
%\subfigure[$h=200, T=1\times 10^5$]{
 \begin{minipage}[t]{0.49\linewidth}
    \centering
     \includegraphics[width=1\textwidth]{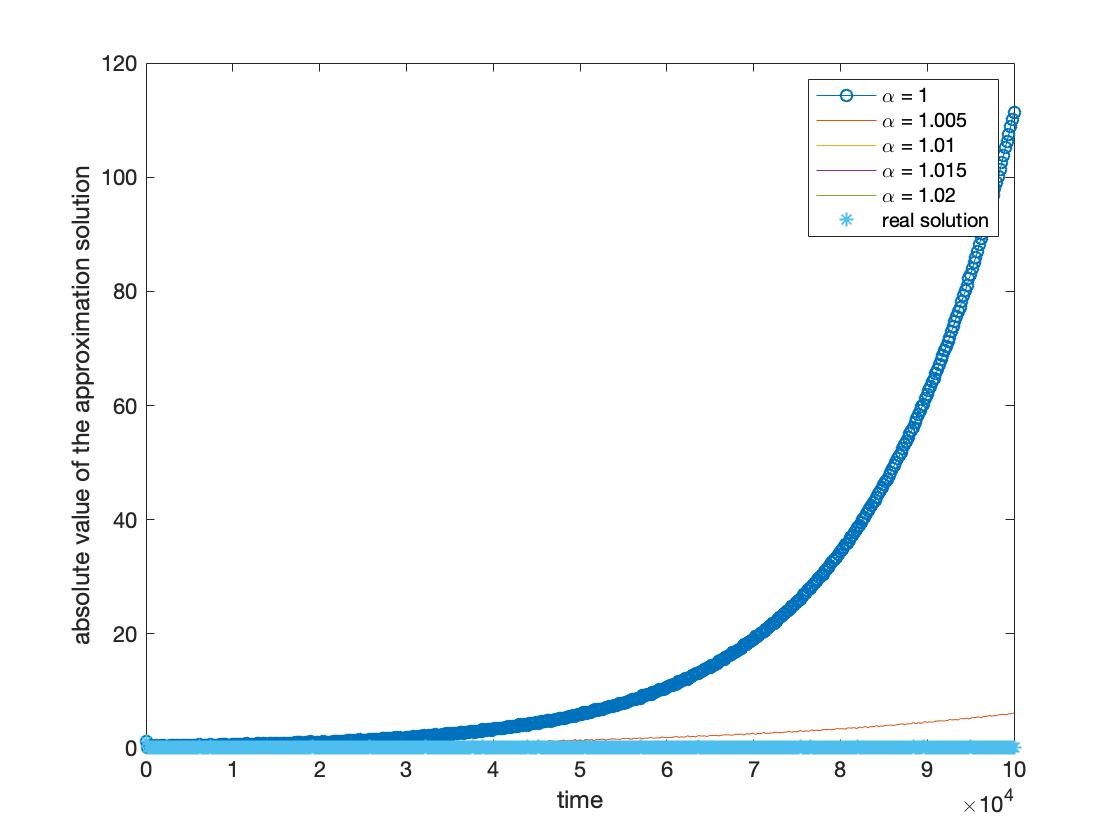}
    %\captionsetup{font={footnotesize}}
%    \caption{$h=200, T=1\times 10^5$}
    \end{minipage}%
%}
\caption{Solution of the generalized BDF2 scheme with different parameter $\alpha$. Left panel $h=100, T=2\times 10^4$.
Right panel $h=200, T=1\times 10^5$.}
\end{figure}

%From figure 5, we observe that the generalized BDF2 scheme is more stable (lower growth rate) for larger parameter $\alpha$. And these results are consistent with our theoretical analysis from the previous section.

\subsubsection{Long time error plot for generalized BDF2 and generalized AM2 methods}

%For long-time behaviour of the generalized BDF2 method and generalized AM2 method, we also give some error plots to show the long-time accuracy and stability.

Here, we demonstrate the long-time second order accuracy of our schemes.

In particular, we applied our IMEX BDF2 scheme \eqref{IMEX} to the the 2D damped-forced ODE system with a skew-symmetric term \eqref{damped_driven_skew}
with $\alpha = 1.1$ or $\alpha = 0.6,$  terminal time $T_{end}=100$, and time-step $h=0.2, 0.1, 0.05$.

\begin{figure}[ht]
\centering
%\subfigure[Generalized BDF2 method with $\alpha = 1.1$]{
   \begin{minipage}[t]{0.49\linewidth}
    \centering
    \includegraphics[width=1\textwidth]{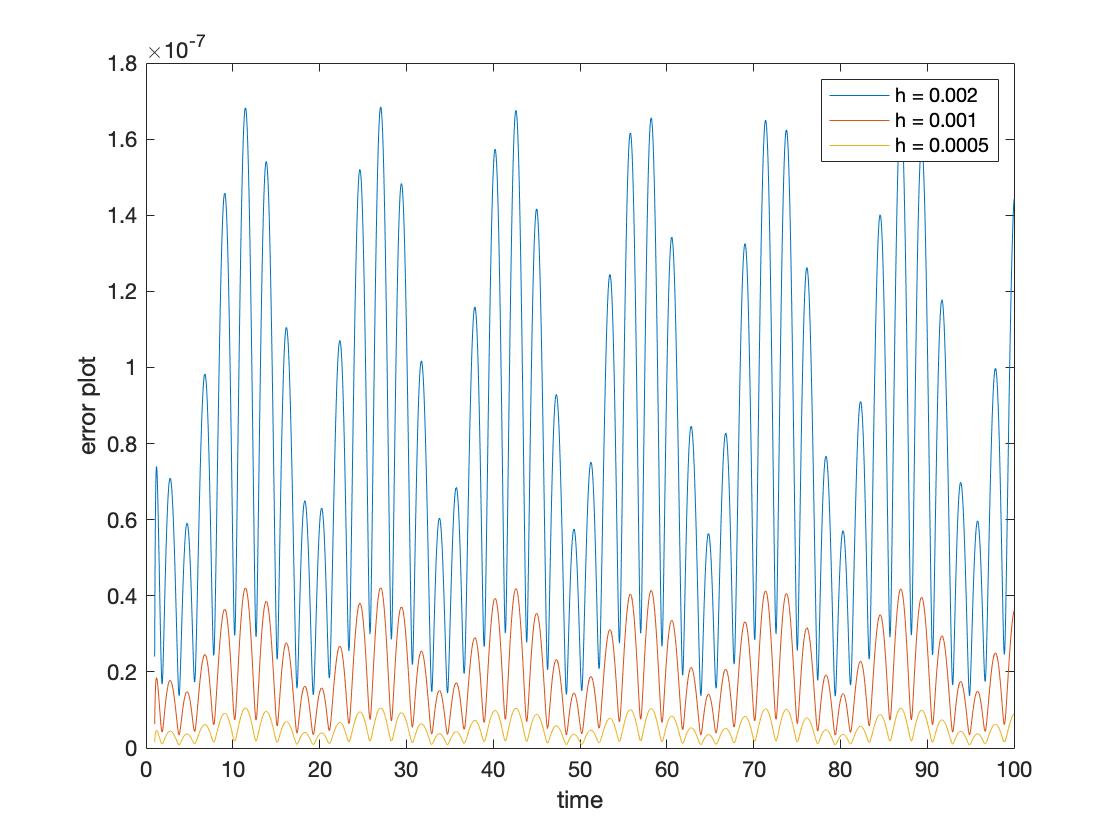}
%    \captionsetup{font={footnotesize}}
%    \caption{Generalized BDF2 method with $\alpha = 1.1$}
    \end{minipage}%
%}
%\subfigure[Generalized AM2 method with $\alpha = 0.6$]{
\begin{minipage}[t]{0.49\linewidth}
    \centering
    \includegraphics[width=1\textwidth]{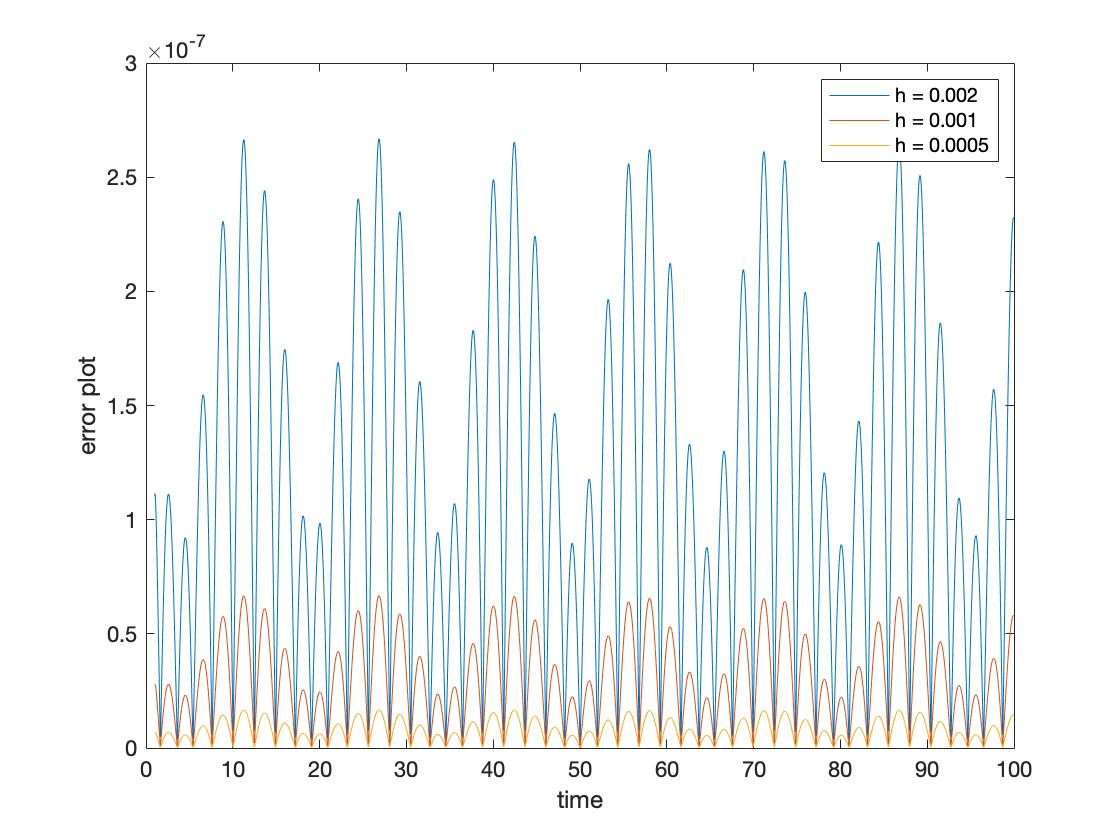}
  %  \captionsetup{font={footnotesize}}
%    \caption{Generalized AM2 method with $\alpha = 0.6$}
    \end{minipage}%
%}
\caption{Relative error over $[0,100]$. Left panel: generalized BDF2-IMEX method with $\alpha = 1.1$. Right panel: Generalized AM2-IMEX method with $\alpha = 0.6$. }
\end{figure}

Tables 5 and 6 demonstrate second order accuracy over a long-time interval (100) for two values of $\alpha$.
Figure 6 contains the relative errors of the numerical solutions for the generalized BDF2 and AM2 method respectively, over the time interval $[0,100]$. We observe that the error is reduced roughly by a factor of 4 when the time step is halved.
In addition, the error is uniformly bounded.

\begin{table}[htbp]
 \caption{long-time error and accuracy for the generalized IMEX BDF2($\alpha = 1.1$) method}
    \centering
    \begin{tabular}{llllll}
\hline
\noalign{\smallskip}
{}&\multicolumn{1}{c}{$h=0.002$} &  \multicolumn{2}{c}{$h= 0.001$}& \multicolumn{2}{c}{$h= 0.0005$}\\
\noalign{\smallskip}
\hline
\noalign{\smallskip}
Evolution Time  &Error &   Error & Order  &Error &  Order\\
\noalign{\smallskip}
\hline
\noalign{\smallskip}
1 & 2.4016e-08  & 6.2603e-09 &  1.9397 & 1.5972e-09 &  1.9707\\

10  &  5.2633e-08 &  1.3079e-08&   2.0087&   3.2600e-09 &  2.0043 \\

20  & 5.8592e-08 &  1.4665e-08   &1.9984&  3.6682e-09 &  1.9992 \\

30 & 1.1102e-07 &  2.7693e-08 &  2.0032 &   6.9153e-09  & 2.0016\\

40 &1.5111e-07 &  3.7803e-08  & 1.9990 &  9.4538e-09  & 1.9995\\

50  &  2.4899e-08&   6.1876e-09  & 2.0086 &   1.5423e-09  & 2.0043\\

60 & 9.6876e-08  & 2.4279e-08 &  1.9964 &  6.0773e-09  & 1.9982\\

70 & 3.7547e-08 &  9.3262e-09  & 2.0093&   2.3240e-09 &  2.0047\\

80 & 4.9994e-08 &  1.2520e-08 &  1.9975&   3.1326e-09  & 1.9988\\

90 & 1.0966e-07 &  2.7344e-08&   2.0038&  6.8270e-09&  2.0019\\

100 &  1.4411e-07 &  3.6043e-08   &1.9994&   9.0126e-09 &  1.9997\\
\noalign{\smallskip}
\hline

\end{tabular}
   
\end{table}

\begin{table}[htbp]
  \caption{long-time error and accuracy for the generalized IMEX AM2($\alpha = 0.6$) method}
    \centering
    \begin{tabular}{llllll}
\hline
\noalign{\smallskip}
{}&\multicolumn{1}{c}{$h=0.002$} &  \multicolumn{2}{c}{$h= 0.001$}& \multicolumn{2}{c}{$h= 0.0005$}\\
\noalign{\smallskip}
\hline
\noalign{\smallskip}
Evolution Time  &Error &   Error & Order  &Error &  Order\\
\noalign{\smallskip}
\hline
\noalign{\smallskip}
1 & 1.1016e-07&   2.7556e-08  & 1.9991&   6.8911e-09 &  1.9996\\

10  &  1.2218e-08 &  2.9711e-09   &2.0399&   7.3246e-10   &2.0202 \\

20  &9.8260e-08 &  2.4568e-08  & 1.9998&   6.1425e-09 &  1.9999\\

30 & 1.3200e-07 &  3.2925e-08   &2.0033&   8.2219e-09  & 2.0016\\

40 &2.4883e-07&   6.2211e-08 &  1.9999&  1.5553e-08  & 2.0000\\

50  &  4.1037e-09 &  1.0003e-09  & 2.0365&  2.4705e-10  & 2.0176\\

60 & 1.8430e-07 &  4.6119e-08&   1.9986&  1.1535e-08 &  1.9993\\

70 & 2.0146e-08 &  5.1253e-09  & 1.9748 &  1.2925e-09 &  1.9875\\

80 & 8.7532e-08 &  2.1892e-08   &1.9994&   5.4740e-09  & 1.9997\\

90 & 1.2336e-07  & 3.0757e-08   &2.0039 &  7.6788e-09  & 2.0020\\

100 & 2.3214e-07 &  5.8029e-08 &  2.0001 &   1.4507e-08  & 2.0001\\
\noalign{\smallskip}
\hline

\end{tabular}
  
\end{table}

% S7
\section{Conclusion and remarks}
\label{Sec7}
We have proposed two families of second order three-level two-step linear schemes that generalize the classical BDF2 and AM2 methods. The generalized schemes are A-stable when the parameter falls into a region that covers the classical cases. 
We have also shown that these schemes are uniform-in-time (energy) stable when applied to simple linear damped-forced models. 
For application to damped-forced model with a low-order skew symmetric term, we proposed IMEX type schemes based on the generalized BDF2/AM2 approach for the time derivative and the dissipation terms while the skew symmetric term is treated explicitly via Gear's extrapolation or Adams-Bashforth approach. 
We are able to show that these generalized BDF2/AM2-IMEX schemes lead to numerical solutions that are uniformly bounded over all time if the skew-symmetric term is relatively weak in some appropriate sense when compared to the linear damping term, and with suitable mild time-step restriction, when the forcing term is uniformly bounded in time.
Our numerics verify our theoretical results. 
Our numerics also indicate that the generalized schemes could be more accurate and more stable than the classical ones. 
%In addition, the stability region of the generalized schemes appears to increase as the parameter increases.
%The rigorous mathematical investigation as well as application to nonlinear problems will be subject to future works.
Note that Dahlquist's second barrier dictates that the highest order of an A-stable multi-step scheme is 2. Therefore, generalization to higher order methods would also be a challenge that calls for different notions of stability.

\section*{Acknowledgements}
This work is based on the thesis of the second author under the supervision of the first author at Southern University of Science and Technology.
The authors acknowledge helpful conversations with Wenbin Chen, Jie Shen, and Jiang Yang. The authors are listed in alphabetic order following convention.
The support  of NSFC and the Havener Endowment is gratefully acknowledged.

%%%%%
\section*{ Appendix: Stokes-Darcy system}
\label{App}

In order to show the potential usage of the generalized BDF2-IMEX as well as the generalized AM2-IMEX schemes, we sketch their application to the 
 physically important Stokes-Darcy system that models flows in fluid saturated karst aquifers in this appendix.
 
 The application to two-dimensional rotation fluid is relatively straightforward in the case of periodic spatial boundary condition.

Recall that the Stokes-Darcy system\cite{cao2010coupled,discacciati2003analysis} can be written as 
\begin{equation}\label{SD-original}
%{\color{red}copy the equation here}
\left\{
\begin{aligned}
    &   S \frac{\partial \phi}{\partial t} - \nabla \cdot (\mathbb{K} \nabla \phi) = f  & \text{ in } \Omega_p\\
&    \frac{\partial \mathbf{u}_f}{\partial t} -\nabla \cdot \mathbb{T}(\mathbf{u}_f,p)  =\mathbf{f} \text{ and } \nabla \cdot \mathbf{u}_f = 0 &\text{ in } \Omega_f,
\end{aligned}
\right.
\end{equation}
with the Beavers-Joseph-Saffman-Jones interface conditions\cite{beavers1967boundary,jones1973low,saffman1971} imposed on the interface $\Gamma$:
\begin{equation}\label{BJSJ_interface_condition}
\left\{    \begin{aligned}
    &    \mathbf{u}_f \cdot \mathbf{n}_f = \mathbf{u}_p \cdot \mathbf{n}_f = -(\mathbb{K}\nabla \phi)\cdot \mathbf{n}_f \\
    & -\mathbf{\tau}_j \cdot (\mathbb{T}(\mathbf{u}_f,p_f)\cdot \mathbf{n}_f) = \alpha_{BJSJ}\mathbf{\tau}_j\cdot \mathbf{u}_f, j=1,2,...,d-1\\
    & -\mathbf{n}_f\cdot (\mathbb{T}(\mathbf{u}_f,p_f)\cdot \mathbf{n}_f) = g\phi,
\end{aligned}
\right.
\end{equation}
where $\mathbf{u}_f$ is the fluid velocity, $p$ is the kinematic pressure, $\phi$ is the hydraulic head, $\mathbf{u}_p = -\mathbb{K}\nabla\phi$ is the velocity in the matrix, $\mathbf{f}$ and $f$ are external body forces acting on the domains $\Omega_f$ and $\Omega_p$, $\mathbb{T}(\mathbf{v},p) = \nu(\nabla\mathbf{v}+\nabla^T\mathbf{v})-p\mathbb{I}$ is the stress tensor, $S$ is the water storage coefficient, $\mathbb{K}$ is the hydraulic conductivity, $\nu$ is the kinematic viscosity of the fluid, $\mathbf{n}_f$ is the outer unit normal vector for $\Omega_f$, $\{ \mathbf{\tau}_j\}_{j=1,2,...,d-1}$ is a linearly independent set of vectors tangent to the interface $\Gamma$, $g$ is the gravitational constant and $\alpha_{BJSJ}$ is the Beavers-Joseph-Saffman-Jones coefficient.

The weak formulation of the Stokes-Darcy system\eqref{SD-original} can be derived by multiplying the equations by test functions $\mathbf{v} \in \mathbf{H}_f, g\phi \in \mathbf{H}_p, q\in Q$ respectively and then integrating over the corresponding domains and then integrating by parts the terms involving second-derivative operators, and then substituting the interface conditions\eqref{BJSJ_interface_condition} in the appropriate terms. The weak formulation\cite{cao2010coupled,chen2013efficient,discacciati2002mathematical} with standard $L^2(D)$ inner product $(\cdot,\cdot)_D$ and norm $\|\cdot\|_D$ and spaces defined by
\begin{equation}
\begin{split}
    \mathbf{H}_f & = \{ \mathbf{v} \in (H^1(\Omega_f))^d | \mathbf{v} = \mathbf{0} \text{ on } \partial \Omega_f \backslash \Gamma \}\\
    H_p & = \{ \psi \in H^1(\Omega_p) | \psi = 0 \text{ on }\partial\Omega_p \backslash \Gamma\} \\
    Q & = L^2(\Omega_f), \quad \mathbf{W} = \mathbf{H}_f \times H_p,\\
\end{split}
\end{equation}
is given as follows: given $f \in (H_p)'$ and $\mathbf{F}\in (\mathbf{H}_f)'$, seek $\phi$ in $H_p, \mathbf{u}_f \in \mathbf{H}_f$, and $p \in Q$ with $\frac{\partial \phi}{\partial t} \in (H_p)'$ and $\frac{\partial \mathbf{u}}{\partial t} \in (\mathbf{H}_f)'$, satisfying
\begin{equation}\label{SD_weak_formulation}
    \begin{split}
        \langle \langle \Vec{\mathbf{u}}_t,\Vec{\mathbf{v}}\rangle \rangle  + a(\Vec{\mathbf{u}},\Vec{\mathbf{v}}) + b(\mathbf{v},q) +a_{\Gamma}(\Vec{\mathbf{u}},\Vec{\mathbf{v}})& = \langle \langle \langle \Vec{\mathbf{f}} , \Vec{\mathbf{v}} \rangle \rangle \rangle\\
        b(\mathbf{u},q) & = 0,
    \end{split}
\end{equation}
where $\Vec{\mathbf{u}} = [\mathbf{u},\phi]^T , \Vec{\mathbf{v}} = [\mathbf{v},\psi]^T$, and $
\Vec{\mathbf{f}} = [\mathbf{f},gf]^T$, $(\cdot)_t = \frac{\partial(\cdot)}{\partial t}$, and notations
\begin{equation}\label{SD_weak_formulation_notation}
    \begin{split}
    &    \langle \langle \Vec{\mathbf{u}}_t,\Vec{\mathbf{v}}\rangle \rangle  = \langle \mathbf{u}_t, \mathbf{v}\rangle_{\Omega_f} + gS\langle \phi_t,\psi \rangle_{\Omega_p}, \quad b(\mathbf{v},q) = \int_{\Omega_f} (-q\cdot \nabla \cdot \mathbf{v})\\
    & a(\Vec{\mathbf{u}},\Vec{\mathbf{v}}) = \nu (\nabla \mathbf{u},\nabla \mathbf{v})_{\Omega_f} +  g(\mathbb{K}\nabla\phi,\nabla\psi)_{\Omega_p} + \alpha_{BJSJ} \int_{\Gamma} \mathbf{u} \cdot \Vec{\tau} \cdot \mathbf{v}\cdot \Vec{\tau}\\
    & a_{\Gamma}(\Vec{\mathbf{u}},\Vec{\mathbf{v}}) =
    g\int_{\Gamma}(\phi \cdot \mathbf{v} \cdot \mathbf{n}_f-\mathbf{u} \cdot \mathbf{n}_f \cdot \psi),\quad  \langle \langle \langle \Vec{\mathbf{f}} , \Vec{\mathbf{v}} \rangle \rangle \rangle = \int_{\Omega_f} \mathbf{f} \cdot \mathbf{v}  + g\int_{\Omega_p} f\psi. \\
    %&a_f(\mathbf{u},\mathbf{v}) = \nu (\nabla \mathbf{u},\nabla \mathbf{v})_{\Omega_f}, a_p(\phi,\psi) = g(\mathbb{K}\nabla\phi,\nabla\psi)_{\Omega_p}, a_{BJSJ}(\mathbf{u},\mathbf{v})= \alpha_{BJ}(\mathbf{u} \cdot \Vec{\tau}\cdot \mathbf{v}\cdot \Vec{\tau})_{\Gamma}.
    \end{split}
\end{equation}

To fit in the structure of the long-time energy stability theorem above, we introduce the following new spaces and notations:

\begin{equation}\label{SD-spaces}
    \begin{aligned}
        H &  = \{ (\bv,\psi)| \bv \in (L^2(\Omega_p))^d, \nabla \cdot \bv = 0, \bv \cdot \bn = 0 \text{ on } \partial \Omega_f \backslash \Gamma, \psi \in L^2(\Omega_p) \}\\
        V & = \{  (\bv,\psi)| \bv \in (H^1(\Omega_p))^d,  \nabla \cdot \bv = 0, \bv  = \mathbf{0} \text{ on } \partial \Omega_f \backslash \Gamma, \psi \in H^1(\Omega_p), \psi = 0 \text{ on } \partial\Omega_p \backslash \Gamma \}
    \end{aligned}
\end{equation}

We denote $V'$ as the dual of $V$ induced by the inner product on $H$:
$$ \left( \begin{pmatrix}
    \bv \\ \phi
\end{pmatrix}, \begin{pmatrix}
    \bu \\ \psi
\end{pmatrix}\right) = \int_{\Omega_f} \bu \cdot \bv + gS\int_{\Omega_p}\phi \psi . $$

\begin{proposition}
     The generalized IMEX BDF2 method \eqref{IMEX-WeakForm} applied on the Stokes-Darcy system \eqref{SD_weak_formulation} on spaces $V$ and $H$ \eqref{SD-spaces} is uniform-in-time energy stable with the mild time step restriction
    \begin{equation}
      \frac{hl_0(\alpha-1)^2}{2}+\frac{h(\alpha-1)^2}{2}+\frac{h\alpha^2}{2C_1} 
\le \frac{4\alpha-3}{4}.
\end{equation}
where $l_0,C_1$ are fixed constants dependent on the system\eqref{SD-original}.
\end{proposition}

    We utilize the stability theorem of the generalized BDF2 and AM2 methods to show the long-time energy stability of the numerical solutions. It is sufficient to show that the assumptions in the stability theorem hold in Stokes-Darcy system. The positive definite term $a(\cdot,\cdot)$ corresponds to a positive definite operator $\mathcal{L}$, s.t. $$a(\Vec{\mathbf{u}},\Vec{\mathbf{v}}) = \langle \mathcal{L}\Vec{\mathbf{u}},\Vec{\mathbf{v}} \rangle,$$ for $\Vec{\bu} ,\Vec{\bv}   \in V$. The verification is followed by utilizing the Poincar\'{e} inequalities
    $$\langle \mathcal{L}y,y \rangle \ge C(\|\bu\|_{H^1}^2 + \|\phi\|_{H^1}^2) \ge l_0\|\begin{pmatrix}
        \bu \\ 
        \phi
    \end{pmatrix} \|_H^2 = l_0\|y\|_H^2,$$
    where $y = \begin{pmatrix}
        \bu \\ 
        \phi
    \end{pmatrix}$.
    
    The skew symmetric term $a_{\Gamma}(\cdot,\cdot) $ corresponds to a skew symmetric operator $\mathcal{L}_s$, s.t.
    $$a_{\Gamma}(\Vec{\mathbf{u}},\Vec{\mathbf{v}}) = \langle \mathcal{L}_s\Vec{\mathbf{u}},\Vec{\mathbf{v}} \rangle,$$
    for $\Vec{\bu} ,\Vec{\bv} \in V$. Notice the skew-symmetry of $\mathcal{L}_s$ as observed in Quarteroni's work. Since $\mathcal{L}_s$ involves the trace of a function, and since the trace is dominated by the $H^1$ norm of the function,and since the weak form of the dissipation term is equivalent to the $H^1$ norm, we have the dominance of the dissipative term verified. Let $y = \begin{pmatrix}
        \bu \\ 
        \phi
    \end{pmatrix}, \tilde{y} = \begin{pmatrix}
        \tilde{\bu} \\ 
        \tilde{\phi}
    \end{pmatrix}$.
    \begin{equation}
        \begin{aligned}
            | \langle \mathcal{L}_sy,\tilde y\rangle|  = & |g\int_{\Gamma}(\phi \cdot \tilde{\bu} \cdot \mathbf{n}_f-\bu \cdot \mathbf{n}_f \cdot \tilde{\phi}) | \\
            \lesssim & \|\phi\|_{H^{\frac{1}{2},00}(\Gamma)} \|\tilde{\bu} \cdot \bn_f\|_{H^{-\frac{1}{2}}(\Gamma)} + \|\tilde{\phi}\|_{H^{\frac{1}{2},00}(\Gamma)} \|\bu \cdot \bn_f\|_{H^{-\frac{1}{2}}(\Gamma)}\\
            \lesssim & \| \phi \|_{H^1(\Omega_p)} \|\tilde{\bu} \|_{L^2(\Omega_f)} + \| \tilde{\phi} \|_{H^1(\Omega_p)} \|\bu \|_{L^2(\Omega_f)}\\
            \lesssim & \|\begin{pmatrix}
        \bu \\ 
        \phi
    \end{pmatrix}\|_V \cdot \|\begin{pmatrix}
        \tilde{\bu} \\ 
        \tilde{\phi}
    \end{pmatrix}\|_H.
        \end{aligned}
    \end{equation}
    
    Therefore we have 
    $$\mathcal{L}_s y \in H, \forall y \in V,$$ and 
    $$\|\mathcal{L}_s y\|_H \lesssim \|y\|_V \lesssim \langle \mathcal{L}y,y \rangle^{\frac{1}{2}}.$$
    
    Equivalently, there exists a constant $C_0>0$, s.t.
    $$\|\mathcal{L}_s y\|_H^2 \le C_0 \langle \mathcal{L}y,y \rangle , \forall y \in V.$$

    Finally, theorem 5 implies the long-time energy stability result for the generalized IMEX BDF2. And similarly, we have the generalized IMEX AM2 case by utilizing theorem 6.
\begin{proposition}
    For $\alpha\in (1/2, 1)$, the generalized IMEX AM2 method \eqref{gAMB2} applied on the Stokes-Darcy system \eqref{SD_weak_formulation} on the corresponding spaces $V$ and $H$ is uniform-in-time energy stable when the mild time step restriction \eqref{gAMB2-step-restriction} holds with parameters depending on the system \eqref{SD_weak_formulation}
\end{proposition}
    
   \noindent \textbf{Remark:} The numerical implementation of the new schemes would naturally involve mixed finite element method. The proof of unconditional stability and uniform-in-time energy bound can be extended to that setting similar to our earlier results \cite{chen2013efficient}.



\end{document}